\documentclass{amsart}

\usepackage{amsmath,amsthm,amssymb}
\usepackage{pb-diagram,lamsarrow,pb-lams}

\author{Qayum Khan}
\title[On fibering and splitting of 5-manifolds over the circle]{On fibering and splitting of \\ 5-manifolds over the circle}
\address{Department of Mathematics, Vanderbilt University, Nashville, TN 37240 U.S.A.}
\email{qayum.khan@vanderbilt.edu}

\newtheorem{thm}{Theorem}[section]

\newtheorem{cor}[thm]{Corollary}

\newtheorem{prop}[thm]{Proposition}
\newtheorem{hyp}[thm]{Hypothesis}

\newtheorem*{unthm}{Theorem}

\theoremstyle{definition}
\newtheorem{defn}[thm]{Definition}

\theoremstyle{remark}
\newtheorem{rem}[thm]{Remark}
\newtheorem*{unrem}{Remark}

\numberwithin{equation}{thm} 
\numberwithin{figure}{section}

\newcommand{\nc}[2]{\newcommand{#1}{#2}}

\nc{\bK}{\mathbf{K}}
\nc{\sCW}{\mathsf{CW}}

\nc{\Cl}{\mathrm{Cl}}
\nc{\Diffeo}{\mathrm{Diffeo}}
\nc{\Homeo}{\mathrm{Homeo}}
\nc{\Hur}{\mathrm{Hur}}
\nc{\Torus}{\mathrm{Torus}}
\nc{\G}{\mathrm{G}}
\nc{\GTOP}{\mathrm{G/TOP}}
\nc{\cS}{\mathcal{S}}
\nc{\ks}{\mathrm{ks}}
\nc{\hofiber}{\mathrm{hofiber}}
\nc{\rel}{\;\mathrm{rel}\,}
\nc{\Hurewicz}{\mathrm{Hurewicz}}
\nc{\cN}{\mathcal{N}}
\nc{\Ker}{\mathrm{Ker}}
\nc{\Cok}{\mathrm{Cok}}
\nc{\Img}{\mathrm{Im}}
\nc{\UNil}{\mathrm{UNil}}
\nc{\CAT}{\mathrm{CAT}}
\nc{\Nil}{\mathrm{Nil}}
\nc{\Wh}{\mathrm{Wh}}
\nc{\Split}{\mathrm{split}}
\nc{\nncc}{\mathrm{nncc}}
\nc{\Id}{\boldsymbol{1}}
\nc{\Isom}{\mathrm{Isom}}
\nc{\PL}{\mathrm{PL}}

\nc{\CP}{\mathbb{CP}}
\nc{\qE}{\mathbb{E}}
\nc{\F}{\mathbb{F}}
\nc{\qH}{\mathbb{H}}
\nc{\R}{\mathbb{R}}
\nc{\RP}{\mathbb{RP}}
\nc{\Z}{\mathbb{Z}}
\nc{\CH}{\mathbb{CH}}
\nc{\bbS}{\mathbb{S}}
\nc{\qL}{\mathbb{L}}
\nc{\act}{\mathrm{act}}
\nc{\hAut}{\mathrm{hAut}}

\nc{\bdry}{\partial}
\nc{\DIFF}{\mathrm{DIFF}}
\nc{\TOP}{\mathrm{TOP}}
\nc{\inv}{^{-1}}
\nc{\homeo}{\approx}
\nc{\x}{\times}
\nc{\en}{\enspace}
\nc{\ST}{\;|\;}
\nc{\longra}{\longrightarrow}
\nc{\xo}{\otimes}
\nc{\iso}{\cong}
\nc{\hookra}{\hookrightarrow}

\newcommand{\wt}[1]{{\widetilde{#1}}}
\newcommand{\xra}[1]{\xrightarrow{#1}}
\newcommand{\gens}[1]{\left\langle #1 \right\rangle}
\newcommand{\SmMatrix}[1]{\left(\begin{smallmatrix} #1\end{smallmatrix}\right)}
\newcommand{\prn}[1]{\left( #1 \right)}
\newcommand{\ol}[1]{\overline{#1}}
\newcommand{\wh}[1]{{\widehat{#1}}}
\newcommand{\set}[1]{{\left\{\, #1 \,\right\}}}

\begin{document}

\begin{abstract}
Our main result is a generalization of Cappell's 5-dimensional
splitting theorem. As an application, we analyze, up to internal
$s$-cobordism, the smoothable splitting and fibering problems for
certain 5-manifolds mapping to the circle. For example, these maps
may have homotopy fibers which are in the class of finite connected
sums of certain geometric 4-manifolds. Most of these homotopy fibers
have non-vanishing second mod 2 homology and have fundamental groups
of exponential growth, which are not known to be tractable by
Freedman--Quinn topological surgery.  Indeed, our key technique is
topological cobordism, which may not be the trace of surgeries.
\end{abstract}

\maketitle 

\section{Introduction}

The problem of whether or not a continuous mapping $f: M \to S^1$ to
the circle from a closed manifold $M$ of dimension $>5$ is homotopic
to a fiber bundle projection was solved originally in the thesis of
F.~Thomas Farrell (cf.~\cite{Farrell_fibering}). The sole
obstruction lies in the Whitehead group of the fundamental group
$\pi_1 M$ and has been reformulated in several ways
\cite{SiebenmannFibering, Farrell_ICM1970, HughesRanicki}. Precious
little is known about the 5-dimensional fibering problem. The
purpose of this paper is to provide more information using recent
advances in rigidity.  Our approach here blends together the
systematic viewpoint of high-dimensional surgery theory and the more
ad-hoc vanishing results known for certain geometric 4-manifolds.

First, we extend some surgery theory.  The central theorem of this
paper is a generalization of the Cappell--Weinberger theorem
\cite{CappellSplit, WeinbergerFibering} for splitting compact
5-manifolds along certain incompressible, two-sided 4-submanifolds
(Theorem \ref{Thm_Split5}). Indeed, the development of additional
tools for our main splitting theorem motivated the author's initial
investigation of 4-manifolds \cite{Khan_Smoothable4}.

Then, we attack the fibering problem. A first application is a
version of the Farrell fibering theorem for smooth $s$-block bundles
(Definition \ref{Defn_sBlock}, Theorem
\ref{Thm_SmoothFiberingOverCircle}) over the circle $S^1$ with
homotopy fiber $\RP^4$ (\ref{Hyp_Projective}); compare
\cite{CS2,FS1,HKT}. The more central geometric applications are to
topological $s$-block bundles (Theorem
\ref{Thm_FiberingOverCircle}). Namely, we allow the fibers to be
compact, orientable 4-manifolds whose interiors admit a complete,
finite volume metrics of  euclidean, real hyperbolic, or complex
hyperbolic type (\ref{Hyp_Surface}). Moreover, we allow the fiber to
be a finite connected sum of orientable surface bundles over
surfaces of positive genus, and of $H$-bundles over the circle $S^1$
such that the compact irreducible 3-manifold $H$ either is $S^3$ or
$D^3$, or is orientable with non-zero first Betti number (hence
Haken), or has complete, finite volume hyperbolic interior
(\ref{Hyp_Sum}).  The hypotheses require smoothness of the total
space and the conclusions assert smoothability of the fiber.

\subsection{Examples of fibers}

Our examples are chosen so that Farrell's fibering obstruction in
$K$-theory and Cappell's splitting obstruction in $L$-theory vanish.

\subsubsection*{First family of examples}
These fibers are certain non-orientable, smooth 4-manifolds with
fundamental group cyclic of order two \cite{HKT}. Assume:

\begin{hyp}\label{Hyp_Projective}
Suppose $Q$ is a non-orientable $\DIFF$ 4-manifold of the form
\[
Q = Q_0 \# Q_1
\]
where:
\begin{enumerate}
\item
$Q_0 = \# r(S^2 \x S^2)$ for some $r \geq 0$, and

\item
$Q_1 = S^2 \x \RP^2$ or $Q_1 = S^2 \rtimes \RP^2$ or $Q_1 =
\#_{S^1} n(\RP^4)$ for some $1 \leq n \leq 4$.
\end{enumerate}
\end{hyp}

\subsubsection*{Second family of examples}

These irreducible, possibly non-orientable fibers have torsion-free
fundamental groups of exponential growth and have non-vanishing
second homotopy group. Assume:

\begin{hyp}\label{Hyp_Surface}
Suppose $S$ is a compact, connected $\DIFF$ 4-manifold such that:
\begin{enumerate}
\item
$S$ is the total space of a $\DIFF$ fiber bundle $S^2 \to S \to
\Sigma$, for some compact, connected, possibly non-orientable
2-manifold $\Sigma$ of positive genus, or

\item
$S$ is the total space of a $\DIFF$ fiber bundle $H \to S \to
S^1$, for some closed, connected, hyperbolic 3-manifold $H$, or

\item
the interior $S-\bdry S$ admits a complete, finite volume metric
of euclidean, real hyperbolic, or complex hyperbolic type.
\end{enumerate}
Moreover, assume $H_1(S;\Z)$ is 2-torsionfree if $S$ is
non-orientable.  Furthermore, in the fiber bundle $S \rtimes_h S^1$
(resp.~$S \rtimes_\alpha S^1$) considered for types (2) and (3) in
this section, assume $h: S \to S$ (resp.~$\alpha$) is homotopic rel
$\bdry S$ to an isometry of $S-\bdry S$.\footnote{This hypothesis is
required since Mostow rigidity fails for product geometries: the
$\qE^4$-manifold $T^2 \x T^2$ has monodromies made from
non-conformal elements of $\pi_0\,\Homeo(T^2)=PSL_2(\Z)$.}
\end{hyp}

\begin{rem}
According to \cite[Lemma 5.9]{Hillman4Book}, the isomorphism classes
of fiber bundles $S^2 \to S \to \Sigma$ in type (1) are in bijective
correspondence with the product $H^1(\Sigma;\Z_2) \x
H^2(\Sigma;\Z_2)$.  The orientable $S^2$-bundles over $\Sigma$ are
classified by the second factor. The isomorphism classes of fiber
bundles $H \to S \to S^1$ in type (2) are in bijective
correspondence with $\pi_0(\Isom\,H)$.
\end{rem}

\subsubsection*{Third family of examples}

These reducible, orientable fibers have torsion-free fundamental
groups of exponential growth and have vanishing second homotopy
groups. A simple example of such a fiber is $F = \# n(S^3 \x S^1)$,
whose fundamental group $\pi_1(F)$ is the free group of rank $n$.
Assume:

\begin{hyp}\label{Hyp_Sum}
Suppose $F$ is an orientable $\DIFF$ 4-manifold of the form
\[
F = F_1 \# \cdots \# F_n
\]
for some $n > 0$, under the following conditions on the compact,
connected, orientable 4-manifolds $F_i$.  Assume:
\begin{enumerate}
\item
$F_i$ is the total space of a $\DIFF$ fiber bundle $H_i \to F_i
\to S^1$, for some compact, connected, orientable 3-manifold
$H_i$ such that:
\begin{enumerate}
\item
$H_i$ is $S^3$ or $D^3$, or

\item
$H_i$ is irreducible with non-zero first Betti number, or
\end{enumerate}

\item
$F_i$ is the total space of a $\DIFF$ fiber bundle $\Sigma_i^f
\to F_i \to \Sigma_i^b$, for some compact, connected, orientable
2-manifolds $\Sigma_i^f$ and $\Sigma_i^b$ of positive genus.
\end{enumerate}
\end{hyp}

\subsection{Main results}

The first splitting theorem is a specialization of the general
splitting theorem (Theorem \ref{Thm_Split5}) to the mapping torus $X
\rtimes_h S^1$ of a homotopy self-equivalence $h: X \to X$ for
certain classes of smooth 4-manifolds $X$.

\begin{unthm}[\ref{Thm_SplittingOverCircle}]
Let $X$ be any of the 4-manifolds $Q,S,F$ defined in
(\ref{Hyp_Projective}, \ref{Hyp_Surface}, \ref{Hyp_Sum}). Let $h: X
\to X$ be a homotopy equivalence which restricts to a diffeomorphism
$\bdry h: \bdry X \to \bdry X$. Suppose $M$ is a compact $\DIFF$
5-manifold and $g: M \to X \rtimes_h S^1$ is a homotopy equivalence
which restricts to a diffeomorphism $\bdry g: \bdry M \to \bdry X
\rtimes_{\bdry h} S^1$.

Then $g$ is homotopic to a map $g'$ which restricts to a simple
homotopy equivalence $g': X' \to X$ such that the $\TOP$ inverse
image $X' := (g')\inv(X)$ is homeomorphic to $X$ and the exterior
$M'$ of $X'$ in $M$ is a smoothable $\TOP$ self $s$-cobordism of
$X$.
\end{unthm}

In other words, cutting $M$ along the bicollared, smoothable $\TOP$
4-submanifold $X' := (g')\inv(X) \homeo X$ yields a smoothable
$\TOP$ $s$-cobordism $(M';X,X)$ and a simple homotopy equivalence
$(g'_\infty;g'_0,g'_1): (M';X,X) \to X \x (\Delta^1;0,1)$ of
manifold triads such that $g'_1 = \alpha \circ g'_0$. Be aware that
the existence of a smooth structure on $X' \homeo X$ does not imply
that $X'$ is a $\DIFF$ submanifold of $M$.

The second splitting theorem connects homotopy structures on mapping
tori to smoothable $s$-cobordisms, homotopy self-equivalences, and
the smoothing invariant.

\begin{unthm}[\ref{Thm_StructuresOverCircle}]
Let $X$ be any of the 4-manifolds $Q,S,F$ defined in
(\ref{Hyp_Projective}, \ref{Hyp_Surface}, \ref{Hyp_Sum}). Let
$\alpha: X \to X$ be a diffeomorphism.  Then there is an exact
sequence of based sets:
\[
\pi_1^\alpha(\wt{\bbS}^s_{\TOP+}(X), \wt{\G}^s(X)) \xra{\en\cup\en}
\cS_\TOP^h(X \rtimes_\alpha S^1) \xra{\en\ks\en} \F_2 \oplus
H_1(X;\F_2)_\alpha.
\]
\end{unthm}

Our fibering theorem is proven using a key strategy of Tom Farrell
\cite{Farrell_ICM1970}. If the smooth 4-manifold $X$ is closed and
simply-connected, the analogous theorem was proven by J.~Shaneson
\cite[Thm.~5.1]{Shaneson_Splitting5}. We do not assume $\bdry X$ is
connected.

\begin{unthm}[\ref{Thm_FiberingOverCircle}]
Let $X$ be any of the 4-manifolds $Q,S,F$ defined in
(\ref{Hyp_Projective}, \ref{Hyp_Surface}, \ref{Hyp_Sum}). Let $M$ be
a $\DIFF$ 5-manifold, and let $f: M \to S^1$ be a continuous map.
Suppose $\bdry X \to \bdry M \xra{\bdry f} S^1$ is a $\DIFF$ fiber
bundle and the homotopy equivalence $\bdry X \to \hofiber(\bdry f)$
extends to a homotopy equivalence $X \to \hofiber(f)$. Then $f: M
\to S^1$ is homotopic $\!\!\rel \bdry M$ to the projection of a
smoothable $\TOP$ $s$-block bundle with fiber $X$.
\end{unthm}

\begin{unrem}
Let $X$ be an aspherical, compact, orientable $\DIFF$ 4-manifold
with fundamental group $\pi$. Suppose the non-connective $L$-theory
assembly map $H_n(\pi;\qL.^h) \to L_n^h(\Z[\pi])$ is an isomorphism
for $n=4,5$. Then the general splitting and fibering theorems
(\ref{Thm_Split5}, \ref{Thm_FiberingOverCircle}) hold for $X$, with
the inclusion of the standard high-dimensional algebraic $K$- and
$L$-theory obstructions.
\end{unrem}

\subsection{Techniques}

Our methods employ geometric topology: topological transversality in
all dimensions (Freedman--Quinn \cite{FQ}) and the prototype of a
nilpotent normal cobordism construction for smooth 5-manifolds
(Cappell \cite{CS1, CappellSplit}). Our hypotheses are
algebraic-topological in nature and come from the surgery
characteristic class formulas of Sullivan--Wall \cite{Wall} and from
the assembly map components of Taylor--Williams \cite{TW}. For the
main application, the difficulty is showing that vanishing of
algebraic $K$- and $L$-theory obstructions is \emph{sufficient} for
a solution to the topological fibering problem as an $s$-block
bundle over the circle.

The reader should be aware that the topological transversality used
in Section \ref{Sec_Assembly5} produces 5-dimensional $\TOP$ normal
bordisms $W\to X \x \Delta^1$ which may not be smoothable, although
$\bdry W = \bdry_- W \cup \bdry_+ W$ is smoothable. In particular,
$W$ may not admit a $\TOP$ handlebody structure relative to $\bdry_-
W$. Hence $W$ may not be the trace of surgeries on topologically
embedded 2-spheres in $X$. Therefore, $W$ may not be produced by
Freedman--Quinn surgery theory, which is developed only for
fundamental groups $\pi_1(X)$ of class $SA$, containing
subexponential growth \cite{KrushkalQuinn}.

\section{Five-dimensional assembly on 4-manifolds}\label{Sec_Assembly5}

Let $(X,\bdry X)$ be a based, compact, connected, $\TOP$ 4-manifold
with fundamental group $\pi = \pi_1(X)$ and orientation character
$\omega = w_1(X): \pi \to \Z^\x$.  Recall, for any $\alpha \in \pi$
and $\beta \in \pi_2(X)$, that there is a Whitehead product
$[\alpha,\beta] \in \pi_2(X)$ which vanishes if and only if the loop
$\alpha$ acts trivially on $\beta$. The \textbf{$\pi$-coinvariants}
are the abelian group quotient $\pi_2(X)_\pi :=
\pi_2(X)/\gens{[\alpha,\beta] \ST \alpha \in \pi, \beta \in
\pi_2(X)}$.

\begin{hyp}\label{Hyp_Assembly5}
Suppose that the following homomorphism is surjective:
\[
\SmMatrix{I_1 & \kappa_3}: H_1(\pi;\Z^\omega) \oplus H_3(\pi;\Z_2)
\longra L_5^h(\Z[\pi]^\omega)
\]
and that the following induced homomorphism is injective:
\[
\Hurewicz: \prn{\pi_2(X) \xo \Z_2}_\pi \longra H_2(X;\Z_2).
\]
\end{hyp}

\begin{thm}\label{Thm_Assembly5}
Assume Hypothesis \ref{Hyp_Assembly5}. Then the following surgery
obstruction map is surjective:
\begin{equation}\label{Eqn_ObstrMap5}
\sigma_*: \cN_\TOP(X \x \Delta^1) \longra L_5^h(\Z[\pi]^\omega).
\end{equation}
\end{thm}

Following Sylvain Cappell's work on the Novikov conjecture, Jonathan
Hillman obtained the same conclusion under different,
group-theoretic hypotheses for a square-root closed graph of certain
class of groups \cite[Lem.~6.9]{Hillman4Book}.

\begin{proof}
There is a commutative diagram
\[
\begin{diagram}
\node{\cN_\TOP(X \x \Delta^1)} \arrow{e,t}{\sigma_*} \arrow{s,tb}{\cap [X]_{\qL^.}}{\iso} \node{L_5^h(\Z[\pi]^\omega).}\\
\node{H_5(X;\GTOP^\omega)} \arrow{e,t}{u_*}
\node{H_5(\pi;\GTOP^\omega)} \arrow{n,t}{A_\pi\!\gens{1}}
\end{diagram}
\]
Since $A_\pi\!\gens{1} = I_1 + \kappa_3$ is surjective and $u_1:
H_1(X;\Z^\omega) \to H_1(\pi;\Z^\omega)$ is an isomorphism, it
suffices to show that $u_3: H_3(X;\Z_2) \to H_3(\pi;\Z_2)$ is
surjective.

Consider the Leray--Serre spectral sequence, with $\pi$-twisted
coefficients, of the fibration $\wt{X} \to X \xra{u} B\pi$, where
$\wt{X}$ is the universal cover of $X$. Then the map $u_3$ is an
edge homomorphism with image subgroup $E^\infty_{3,0}$.  Note
\[
E^\infty_{3,0} = \Ker\prn{d^3_{3,0}: H_3(B\pi;\Z_2) \longra
(\pi_2(X) \xo \Z_2)_\pi}.
\]
There is an exact sequence involving the associated graded groups
$E^\infty_{0,2}$ and $E^\infty_{2,0}$ and inducing the classical
Hopf sequence:
\[
0 \xra{\en} \Cok(d^3_{3,0}) \xra{\Hurewicz_*} H_2(X;\Z_2) \xra{u_*}
H_2(B\pi;\Z_2) \xra{\en} 0.
\]
It follows from the second part of the hypothesis that the
transgression $\bdry = d^3_{3,0}$ is zero. Therefore $\Img(u_3) =
E^\infty_{3,0} = H_3(\pi;\Z_2)$, hence $\sigma_*$ is surjective.
\end{proof}

Some families of reducible examples $X$ of the theorem are obtained
as finite connected sums of certain compact, aspherical 4-manifolds
$X_i$ which are constructed from non-positively curved manifolds.
Recall that the interior of any compact surface $\Sigma$ of positive
genus has the structure of a complete, finite volume, euclidean or
hyperbolic 2-manifold; hence $\Sigma$ is aspherical. The following
corollary gives a rich source of examples, including $X = \# n(S^3
\x S^1)$, whose fundamental group is free.

\begin{cor}\label{Cor_Assembly5}
Suppose $X$ is a $\TOP$ 4-manifold of the form
\[
X = X_1 \# \cdots \# X_n
\]
for some $n > 0$ and some compact, connected 4-manifolds $X_i$ with
(torsion-free) fundamental groups $\Lambda_i$ such that:
\begin{enumerate}
\item
the interior $X_i - \bdry X_i$ admits a complete, finite volume
metric of real or complex hyperbolic type, or

\item
$X_i$ is the total space of a fiber bundle $\Sigma_i^f \to X_i
\to \Sigma_i^b$, for some compact, connected 2-manifolds
$\Sigma_i^f$ and $\Sigma_i^b$ of positive
genus\footnote{\textbf{Positive genus:} implies torsion-free;
each surface $\Sigma_i^f$ and $\Sigma_i^b$ is a finite connected
sum of at least either one torus $T^2$ or two real projective
planes $\RP^2$, with arbitrary punctures. The first
non-orientable example is the Klein bottle $Kl = \RP^2 \#
\RP^2$.} , or

\item
$X_i$ is the total space of a fiber bundle $H_i \to X_i \to S^1$,
for some compact, connected, irreducible 3-manifold $H_i$ such that:
\begin{enumerate}
\item
$H_i$ is $S^3$ or $D^3$,

\item
$H_i$ is orientable with non-zero first Betti number
(i.e.~$H^1(H_i;\Z) \neq 0$, e.g.~the boundary $\bdry H_i$ is
non-empty), or

\item
the interior $H_i - \bdry H_i$ admits a complete, finite
volume metric of hyperbolic type.
\end{enumerate}
\end{enumerate}
Then the topological 5-dimensional surgery obstruction map
(\ref{Eqn_ObstrMap5}) is surjective.  Moreover, $X_i$ only needs to
have type (1), (2), or (3) up to homotopy equivalences which respect
orientation characters.
\end{cor}

\begin{proof}
Let $\Lambda_i := \pi_1(X_i)$ be the fundamental group of $X_i$ with
orientation character $\omega_i: \Lambda_i \to \Z^\x$.  Consider the
connective assembly map
\[
A_{\Lambda_i}\!\gens{1} = \SmMatrix{I_1 & \kappa_3}:
H_5(B\Lambda_i;\GTOP^{\omega_i}) \longra
L_5^h(\Z[\Lambda_i]^{\omega_i}).
\]
In order to verify Hypothesis \ref{Hyp_Assembly5} and apply Theorem
\ref{Thm_Assembly5}, it suffices to show that:
\begin{enumerate}
\item[(i)]
$\pi_2(X_i)=0$,

\item[(ii)]
$H_d(B\Lambda_i;\Z)=0$ for all $d> 4$, and

\item[(iii)]
the non-connective assembly map is an isomorphism:
\[
A_{\Lambda_i}: H_5(B\Lambda_i; \qL.^{\omega_i}) \longra
L_5^h(\Z[\Lambda_i]^{\omega_i}).
\]

\end{enumerate}
Then $A_{\Lambda_i}\!\gens{1}$ is an isomorphism. So, since the
trivial group $1$ is square-root closed in the torsion-free groups
$\Lambda_i$, the $\UNil$-groups associated to the free product $\pi
= \bigstar_{i=1}^n \Lambda_i$ vanish, by \cite[Corollary
4]{CappellUnitary}, which was proven in \cite[Lemmas
II.7,8,9]{CappellSplit}.  Therefore, by the Mayer--Vietoris sequence
in $L$-theory \cite[Thm.~5(ii)]{CappellUnitary}, Proposition
\ref{Prop_VanishingWhiteheadGroup} for $h$-decorations, the
five-lemma, and induction on $n$, we obtain that $A_\pi\!\gens{1}$
is an isomorphism. Moreover, by the Mayer--Vietoris sequence in
singular homology and the Hurewicz theorem applied to the universal
cover $\wt{X}$, we obtain that $\pi_2(X)=0$.

There are three types of connected summands $X_i$.

\emph{Type 1.} Since $X_i - \bdry X_i$ is covered by $\qH^4$ or
$\CH^2$, we obtain $X_i$ is aspherical. That is, the compact
4-manifold $X_i$ is model for $B\Lambda_i$. Then (i) and (ii) are
satisfied. Since $X - \bdry X_i$ is complete, homogeneous, and has
non-positive sectional curvatures, by \cite[Proposition
0.10]{FJ_TorsionFreeGL}, condition (iii) is satisfied.

\emph{Type 2.} Since the surfaces $\Sigma_i^f$ and $\Sigma_i^b$ are
aspherical, by the homotopy fibration sequence, we obtain that the
compact 4-manifold $X_i$ is aspherical.  Then (i) and (ii) are
satisfied. By a result of J.~Hillman \cite[Lemma
6]{Hillman_SurfaceBundles} for closed, aspherical surface bundles
over surfaces, condition (iii) is satisfied. Indeed, the
Mayer--Vietoris argument extends to compact, aspherical surfaces
with boundary: each circle $C_j$ in the connected-decomposition of
the aspherical surface $\Sigma_i^b = F_1 \# \cdots \# F_r$ generates
an indivisible element in the free fundamental group of the
many-punctured torus or Klein bottle $F_k$, hence each inclusion
$\pi_1(C_j) \to \pi_1(F_k)$ of fundamental groups is square-root
closed (see \cite[Thm.~2.4]{CHS} for detail).

\emph{Type 3.} There are three types of fibers $H_i$.

\emph{Type 3a.} Conditions (i)--(iii) are immediately satisfied.

\emph{Type 3b.} Since $H_i$ is a compact, connected, irreducible,
orientable 3-manifold and $\pi_1(H_i)$ is infinite, using the Sphere
Theorem of Papakyriakopoulos and the Hurewicz theorem, it can be
shown that $H_i$ is aspherical. Then $X_i$ is aspherical, so
conditions (i) and (ii) are satisfied. Since $H_i$ is irreducible
and $H_i \neq D^3$, no connected component of $\bdry H_i$ is a
2-sphere.  If $\bdry H_i$ is non-empty, then it can be shown that
$H_i$ is Haken \cite[Lem.~6.8]{Hempel3Book}. So, by theorems of S.
Roushon, the non-connective assembly map $A_{\pi_1(H_i)}$ is an
isomorphism in dimensions 4 and 5: if $\bdry H$ is non-empty, this
follows from \cite[Theorem 1.1(1)]{Roushon_Haken}, and if $\bdry H$
is empty, this follows from \cite[Theorem 1.2]{Roushon_Betti}.
Therefore, by the Ranicki--Shaneson sequence
\cite[Thm.~5.2]{RanickiLIII}, Proposition
\ref{Prop_VanishingWhiteheadGroup} for $h$-decorations, and the
five-lemma, we obtain that condition (iii) is satisfied.

\emph{Type 3c.}  Since $\qH^3$ is the universal cover of $H_i -
\bdry H_i$, the interior $H_i - \bdry H_i$ is aspherical. But, since
$\bdry H_i$ has a collar implies that $H_i - \bdry H_i \hookra H_i$
is a homotopy equivalence, we obtain that $H_i$ is also aspherical.
Then $X_i$ is aspherical, so conditions (i) and (ii) are satisfied.
Since $\qH^3$ isometrically covers $H_i - \bdry H_i$, by a result of
Farrell and Jones \cite[Prop.~0.10]{FJ_TorsionFreeGL}, the
non-connective assembly map $A_{\pi_1(H_i)}$ is an isomorphism in
dimensions 4 and 5. Therefore, by the Ranicki--Shaneson sequence
\cite[Thm.~5.2]{RanickiLIII}, Proposition
\ref{Prop_VanishingWhiteheadGroup} for $h$-decorations, and the
five-lemma, we obtain that condition (iii) is satisfied.
\end{proof}

Here is a family of non-aspherical examples $X$ of the theorem.

\begin{cor}\label{Cor_Assembly5Nonaspherical}
Suppose $X$ is a compact $\TOP$ 4-manifold which is homotopy
equivalent to the total space of a fiber bundle $S^2 \to E \to
\Sigma$, for some compact, connected 2-manifold $\Sigma$ of positive
genus.  Then the topological 5-dimensional surgery obstruction map
(\ref{Eqn_ObstrMap5}) is surjective.
\end{cor}

\begin{proof}
By \cite[Theorem 6.16]{Hillman4Book}, $X$ is $s$-cobordant to $E$.
Hence there is an induced (simple) homotopy equivalence $X \to E$
which respects orientation characters. The same methods of Corollary
\ref{Cor_Assembly5}(2) show that $\SmMatrix{I_1 & \kappa_3}$ is
surjective.   Note that $\pi_1(E) = \pi_1(\Sigma)$ may not act
trivially on $\pi_2(E) = \pi_2(S^2) = \Z$ but does acts trivially on
$\pi_2(E) \xo \Z_2 = \Z_2$.  An elementary argument with the
Leray--Serre spectral sequence shows that $H_2(E;\Z_2) =
H_2(S^2;\Z_2) \oplus H_2(\Sigma;\Z_2)$.  Therefore $\Hurewicz$ is
injective, and we are done by Theorem \ref{Thm_Assembly5}.
\end{proof}

The difference between $\DIFF$ and $\TOP$ for $\sigma_*$ is
displayed in \cite[Prop.~2.1]{Khan_Smoothable4}. Later, we shall
refer to a hypothesis introduced by Cappell \cite[Thm.~5,
Rmk.]{CappellSplit}.

\begin{hyp}\label{Hyp_Surj5}
Suppose that the following map is surjective:
\[
\sigma_*: \cN_\DIFF(X \x \Delta^1) \longra L_5^h(\Z[\pi]^\omega).
\]
\end{hyp}

\begin{rem}\label{Rem_OrderTwo}
Suppose $X$ is a $\DIFF$ 4-manifold and $\pi^\omega = (C_2)^-$. By
\cite[Theorem 13A.1]{Wall}, the following surgery obstruction map is
automatically surjective:
\[
\sigma_*: \cN_\DIFF(X \x \Delta^1) \longra L_5^h(\Z[C_2]^-) = 0.
\]
Topological surjectivity fails for a connected sum $X \# X'$ of such
manifolds: the Mayer--Vietoris sequence \cite{CappellUnitary} shows
that the cokernel is $\UNil_5^h(\Z;\Z^-,\Z^-) \iso
\UNil_3^h(\Z;\Z,\Z)$, and this abelian group was shown to be
infinitely generated \cite{CR}.
\end{rem}

\section{Exactness at 4-dimensional normal
invariants}\label{Sec_Exactness}

For the convenience of the reader, we first recall the relevant
hypotheses from the precursor \cite[\S 3]{Khan_Smoothable4}. Let
$(X,\bdry X)$ be a based, compact, connected $\TOP$ 4-manifold with
fundamental group $\pi$ and orientation character $\omega: \pi \to
\Z^\x$. Let $u: X \to B\pi$ be a based, continuous map that induces
an isomorphism on fundamental groups. Denote the induced
homomorphism
\[
u_2: H_2(X;\Z_2) \longra H_2(B\pi;\Z_2).
\]
Recall that $X$ satisfies Poincar\'e duality with a unique mod 2
orientation class $[X] \in H_4(X,\bdry X;\Z_2)$. The second Wu class
$v_2(X) \in H^2(X;\Z_2)$ is a homomorphism
\[
v_2(X): H_2(X;\Z_2) \longra \Z_2
\]
uniquely determined, for all cohomology classes $a \in H^2(X,\bdry
X;\Z_2)$, by the formula
\[
\gens{v_2(X), a \cap [X]} = \gens{a \cup a, [X]}.
\]

Consider three cases for the orientation character $\omega$ below.
The homomorphism
\[
\kappa_2: H_2(B\pi;\Z_2) \longra L_4^h(\Z[\pi]^\omega)
\]
is the 2-dimensional component of the $L$-theory assembly map.

\begin{hyp}\label{Hyp_Orient} Let $X$ be orientable.
Suppose that $\kappa_2$ is injective on the subgroup
$u_2(\Ker\,v_2(X))$.
\end{hyp}

\begin{hyp}\label{Hyp_NonorientFin}
Let $X$ be non-orientable such that $\pi$ contains an
orientation-reversing element of finite order, and if $\CAT=\DIFF$,
then suppose that orientation-reversing element has order two.
Suppose that $\kappa_2$ is injective on all $H_2(B\pi;\Z_2)$, and
suppose that $\Ker(u_2) \subseteq \Ker(v_2)$.
\end{hyp}

\begin{hyp}\label{Hyp_NonorientInf}
Let $X$ be non-orientable such that there exists an epimorphism
$\pi^\omega \to \Z^-$. Suppose that $\kappa_2$ is injective on the
subgroup $u_2(\Ker\,v_2(X))$.
\end{hyp}

Next, we recall the relevant results from \cite[\S
4]{Khan_Smoothable4} used frequently in the later proofs in this
paper. The subcategory $\TOP0 \subset \TOP$ consists of those maps
$f: M \to X$ with Kirby--Siebenmann stable smoothing invariant
$\ks(f) := \ks(M) - \ks(X) =0 \in \Z_2$. All structure sets and
normal invariants below are relative to a diffeomorphism on $\bdry
X$.

\begin{thm}\label{Thm_SES}
Let $(X,\bdry X)$ be a based, compact, connected, $\CAT$ 4--manifold
with fundamental group $\pi = \pi_1(X)$ and orientation character
$\omega = w_1(X): \pi \to \Z^\x$.
\begin{enumerate}
\item
Suppose Hypothesis \ref{Hyp_Orient} or \ref{Hyp_NonorientFin}.
Then the surgery sequence of based sets is exact at the smooth
normal invariants:
\begin{equation}\label{SES_DIFF}
\cS^s_\DIFF(X) \xra{\en\eta\en} \cN_\DIFF(X) \xra{\en\sigma_*\en}
L_4^h(\Z[\pi]^\omega).
\end{equation}
\item
Suppose Hypothesis \ref{Hyp_Orient} or \ref{Hyp_NonorientFin} or
\ref{Hyp_NonorientInf}.  Then the surgery sequence of based sets
is exact at the stably smoothable normal invariants:
\begin{equation}\label{SES_TOP}
\cS^s_{\TOP 0}(X) \xra{\en\eta\en} \cN_{\TOP 0}(X)
\xra{\en\sigma_*\en} L_4^h(\Z[\pi]^\omega).
\end{equation}
\end{enumerate}\qed
\end{thm}

\begin{cor}\label{Cor_FreeProdCrys}
Let $\pi$ be a free product of groups of the form
\[
\pi = \bigstar_{i=1}^n \Lambda_i
\]
for some $n > 0$, where each $\Lambda_i$ is a torsion-free lattice
in either $\Isom(\qE^{m_i})$ or $\Isom(\qH^{m_i})$ or
$\Isom(\CH^{m_i})$ for some $m_i>0$. Suppose the orientation
character $\omega$ is trivial. Then the surgery sequences
(\ref{SES_DIFF}) and (\ref{SES_TOP}) are exact.\qed
\end{cor}

\begin{cor}\label{Cor_OrderTwo}
Suppose $X$ is a $\DIFF$ 4--manifold of the form
\[
X = X_1 \# \cdots \# X_n \# r(S^2 \x S^2)
\]
for some $n > 0$ and $r \geq 0$, and each summand $X_i$ is either
$S^2 \x \RP^2$ or $S^2 \rtimes \RP^2$ or $\#_{S^1} n (\RP^4)$ for
some $1 \leq n \leq 4$. Then the surgery sequences (\ref{SES_DIFF})
and (\ref{SES_TOP}) are exact.\qed
\end{cor}

\begin{cor}\label{Cor_Haken}
Suppose $X$ is a $\TOP$ 4--manifold of the form
\[
X = X_1 \# \cdots \# X_n \# r (S^2 \x S^2)
\]
for some $n>0$ and $r \geq 0$, and each summand $X_i$ is the total
space of a fiber bundle
\[
H_i \longra X_i \longra S^1.
\]
Here, we suppose $H_i$ is a compact, connected 3--manifold such
that:
\begin{enumerate}
\item
$H_i$ is $S^3$ or $D^3$, or

\item
$H_i$ is irreducible with non-zero first Betti number.
\end{enumerate}
Moreover, if $H_i$ is non-orientable, we assume that the quotient
group $H_1(H_i;\Z)_{(\alpha_i)_*}$ of coinvariants is 2-torsionfree,
where $\alpha_i: H_i \to H_i$ is the monodromy homeomorphism. Then
the surgery sequence (\ref{SES_TOP}) is exact.\qed
\end{cor}

\begin{cor}\label{Cor_Klein}
Suppose $X$ is a $\TOP$ 4--manifold of the form
\[
X = X_1 \# \cdots \# X_n \# r(S^2 \x S^2)
\]
for some $n > 0$ and $r \geq 0$, and each summand $X_i$ is the total
space of a fiber bundle
\[
\Sigma_i^f \longra X_i \longra \Sigma_i^b.
\]
Here, we suppose the fiber and base are compact, connected
2--manifolds, $\Sigma_i^f \neq \RP^2$, and $\Sigma_i^b$ has positive
genus. Moreover, if $X_i$ is non-orientable, we assume that the
fiber $\Sigma_i^f$ is orientable and that the monodromy action of
$\pi_1(\Sigma_i^b)$ of the base preserves any orientation on the
fiber. Then the surgery sequence (\ref{SES_TOP}) is exact.\qed
\end{cor}

\section{Splitting of 5-manifolds}

We generalize Cappell's 5-dimensional splitting theorem
\cite[Thm.~5, Remark]{CappellSplit}, using the homological
hypotheses developed in Sections
\ref{Sec_Assembly5}---\ref{Sec_Exactness}. Our proof incorporates
the possible non-vanishing of $\UNil_6$.  The $\DIFF$ and $\TOP$
cases are distinguished, and the results of this section are applied
to the fibering problem in Section \ref{Sec_SplitFiber}. The stable
surgery version of the splitting theorem can be found in \cite{CS1}.
However, the stable splitting of 5-manifolds is not pursued here,
since connecting sum a single fiber with $S^2 \x S^2$ destroys the
fibering property over $S^1$.

Let $(Y,\bdry Y)$ be a based, compact, connected $\CAT$ 5-manifold.
Let $(Y_0,\bdry Y_0)$ is a based, compact, connected $\CAT$
4-manifold. Suppose $Y_0$ is an \textbf{incompressible, two-sided}
submanifold of $Y$. That is, the induced homomorphism $\pi_1(Y_0)
\to \pi_1(Y)$ is injective, and there is a separating decomposition
\[
Y = Y_- \cup_{Y_0} Y_+ \quad\text{with}\quad \bdry Y = \bdry Y_-
\cup_{\bdry Y_0} \bdry Y_+
\]
or, respectively, a non-separating decomposition
\[
Y = \cup_{Y_0} Y_\infty \quad\text{with}\quad \bdry Y = \cup_{\bdry
Y_0}\, \bdry Y_\infty.
\]
The Seifert--van Kampen theorem identifies
\[
\pi_1(Y) = \Pi = \Pi_- *_{\Pi_0} \Pi_+
\]
as the corresponding injective, amalgamated free product of
fundamental groups, or, respectively,
\[
\pi_1(Y) = \Pi = *_{\Pi_0}\, \Pi_\infty
\]
as the corresponding injective, HNN-extension\footnote{In the
non-separating case, we write $\Pi_0^-, \Pi_0^+$ as the two
monomorphic images of $\Pi_0$ in $\Pi_\infty$.} of fundamental
groups.

A homotopy equivalence $g$ to $Y$ is \textbf{$\CAT$ splittable along
$Y_0$} if $g$ is homotopic, relative to a $\CAT$ isomorphism $\bdry
g$, to a union $g_- \cup_{g_0} g_+$ (resp.~$\cup_{g_0}\, g_\infty$)
of homotopy equivalences from compact $\CAT$ manifolds to $Y_-, Y_0,
Y_+$ (resp.~$Y_0, Y_\infty$) \cite{CappellFree}. Under certain
conditions, we show that the vanishing of high-dimensional
obstructions in $\Nil_0$ and $\UNil_6^s$ are sufficient for
splitting. These two obstructions were formulated by Friedhelm
Waldhausen (1960's) and Sylvain Cappell (1970's).

\begin{thm}\label{Thm_Split5}
Let $(Y,\bdry Y)$ be a finite, simple Poincar\'e pair of formal
dimension 5 \cite[\S 2]{Wall}. Suppose $\bdry Y$ and $Y_0$ are
compact $\DIFF$ 4-manifolds such that $(Y_0,\bdry Y_0)$ is a
connected, incompressible, two-sided Poincar\'e subpair of $(Y,\bdry
Y)$ with tubular neighborhood $Y_0 \x [-1,1]$. If $\CAT=\DIFF$,
assume $Y_0$ satisfies Hypothesis \ref{Hyp_Orient} or
\ref{Hyp_NonorientFin} and satisfies Hypothesis \ref{Hyp_Surj5}. If
$\CAT=\TOP$, assume $Y_0$ satisfies Hypothesis \ref{Hyp_Orient} or
\ref{Hyp_NonorientFin} or \ref{Hyp_NonorientInf} and satisfies
Hypothesis \ref{Hyp_Assembly5}.

Suppose $g: (W,\bdry W) \to (Y,\bdry Y)$ is a homotopy equivalence
for some compact $\DIFF$ 5-manifold $W$ such that the restriction
$\bdry g: \bdry W \to \bdry Y$ is a diffeomorphism. Then $g$ is
$\CAT$ splittable along $Y_0$ if and only if
\begin{enumerate}
\item
the cellular splitting obstruction, given by the image of the
Whitehead torsion $\tau(g) \in \Wh_1(\Pi)$, vanishes:
\[
\Split_K(g;Y_0) \in \Wh_0(\Pi_0) \oplus
\wt{\Nil}_0(\Z[\Pi_0];\Z[\Pi_- - \Pi_0], \Z[\Pi_+ - \Pi_0])
\]
or, respectively,
\begin{multline*}
\qquad\quad\en \Split_K(g;Y_0) \in \Wh_0(\Pi_0) \en\oplus\\
\wt{\Nil}_0(\Z[\Pi_0];\Z[\Pi_\infty - \Pi_0^-], \Z[\Pi_\infty -
\Pi_0^+], {}_{-}\Z[\Pi_\infty]_{+}, {}_{+}\Z[\Pi_\infty]_{-})
\end{multline*}
and subsequently
\item
the manifold splitting obstruction, given by the algebraic position
of discs in the fundamental subdomains of the $\Pi_0$-cover,
vanishes:
\[
\Split_L(g;Y_0) \in \UNil_6^s(\Z[\Pi_0]^{\omega_0};\Z[\Pi_- -
\Pi_0]^{\omega_-}, \Z[\Pi_+ - \Pi_0]^{\omega_+})
\]
or, respectively,
\[
\Split_L(g;Y_0) \in \UNil_6^s(\Z[\Pi_0]^{\omega_0};\Z[\Pi_\infty -
\Pi_0^-]^{\omega_\infty}, \Z[\Pi_\infty - \Pi_0^+]^{\omega_\infty}).
\]
\end{enumerate}
Furthermore, if $g$ is $\CAT$ splittable along $Y_0$, then $g$ is
homotopic $\!\!\rel \bdry W$ to a split homotopy equivalence $g': W
\to Y$ such that the $\CAT$ inverse image $(g')\inv(Y_0)$ is $\CAT$
isomorphic to $Y_0$.
\end{thm}

Our theorem mildly generalizes \cite[Theorem 5,
Remark]{CappellSplit}, which included: if $\Pi_0$ is a finite group
of odd order, then $H_2(\Pi_0;\Z_2)=0$ and $L_5^h(\Z[\Pi_0])=0$.

\begin{cor}[Cappell]
Suppose $g: W \to Y$ is a homotopy equivalence of closed $\DIFF$
5-manifolds. Assume:
\begin{enumerate}
\item
$Y_0$ is orientable,
\item
$H_2(\Pi_0;\Z_2)=0$,
\item
$\Pi_0$ is square-root closed\footnote{\textbf{Square-root closed}:
if $g \in \Pi$, then $g^2 \in \Pi_0$ implies $g \in \Pi_0.$} in
$\Pi$, and
\item
the following surgery obstruction map is surjective
(cf.~Hyp.~\ref{Hyp_Surj5}):
\[
\sigma_*: \cN_\DIFF(X \x \Delta^1) \longra L_5^h(\Z[\Pi_0]).
\]
\end{enumerate}

Then $g$ is $\DIFF$ splittable along $Y_0$ if and only if the above
image $\Split_K(g;Y_0)$ of the Whitehead torsion $\tau(g) \in
\Wh_1(\Pi)$ vanishes. \qed
\end{cor}

Define a decoration subgroup $B \subseteq \Wh_1(\Pi)$ as the image
of $\Wh_1(\Pi_-) \oplus \Wh_1(\Pi_+)$, respectively
$\Wh_1(\Pi_\infty)$, under the homomorphism induced by inclusion.
Recall that the \textbf{structure set} $\cS_\CAT^B(Y)$ is defined as
the set of equivalence classes of homotopy equivalences $g: (W,\bdry
W) \to (Y,\bdry Y)$ such that $W$ is a compact $\CAT$ manifold and
$\bdry g: \bdry W \to \bdry Y$ is a $\CAT$ isomorphism and $g$ has
Whitehead torsion $\tau(g) \in B$, under the equivalence relation $g
\sim g'$ if there exists a $\CAT$ isomorphism $h: W \to W'$ such
that $g' \circ h$ is homotopic to $g$. The \textbf{split structure
set} $\cS_\CAT^\Split(Y;Y_0)$ is defined as the subset of
$\cS_\CAT^B(Y)$ whose elements are represented by homotopy
equivalences $\CAT$ splittable along $Y_0$. The abelian group
$\UNil_6^s$ depends only on the fundamental groups $\Pi_-, \Pi_0,
\Pi_+$ (resp.~$\Pi_0,\Pi_\infty$) with orientation character
$\omega$. $\UNil_6^s$ is algebraically defined and has zero
decoration in $\wt{\Nil}_0$ \cite{CappellUnitary}.

\begin{defn}\label{Defn_TOPplus}
Let $(Y,\bdry Y)$ be a compact $\DIFF$ manifold. Define the
\textbf{smoothable structure set} $\cS_{\TOP +}(Y)$ as the image of
$\cS_\DIFF(Y)$ under the forgetful map to $\cS_\TOP(Y)$.  That is,
$\cS_{\TOP+}(Y)$ is the subset of $\cS_\TOP(Y)$ consisting of the
elements representable by homotopy equivalences $g: (W, \bdry W) \to
(Y, \bdry Y)$ such that $W$ admits a $\DIFF$ structure extending the
$\DIFF$ structure on $\bdry W$ induced by $\bdry g$.
\end{defn}

A more succinct statement illuminates the method of proof in higher
dimensions: Sylvain Cappell's ``nilpotent normal cobordism
construction'' \cite{CappellFree, CappellSplit}.

\begin{thm}\label{Thm_nncc}
Let $(Y,\bdry Y)$ be a finite, simple Poincar\'e pair of formal
dimension 5 \cite[\S 2]{Wall}. Suppose $\bdry Y$ and $Y_0$ are
compact $\DIFF$ 4-manifolds such that $(Y_0,\bdry Y_0)$ is a
connected, incompressible, two-sided Poincar\'e subpair of $(Y,\bdry
Y)$ with tubular neighborhood $Y_0 \x [-1,1]$.

\begin{enumerate}
\item
Assume $Y_0$ satisfies Hypothesis \ref{Hyp_Orient} or
\ref{Hyp_NonorientFin} and satisfies Hypothesis \ref{Hyp_Surj5}.
Then there is a bijection
\[
\nncc^s: \cS_\DIFF^B(Y) \longra \cS_\DIFF^\Split(Y;Y_0) \x \UNil_6^s
\]
such that composition with projection onto the first factor is a
subset retraction, and composition with projection onto the second
factor is the manifold splitting obstruction $\Split_L$.
Furthermore, $g$ and $\nncc^s(g)$ have equal image in
$\cN_\DIFF(Y)$.

\item
Assume $Y_0$ satisfies Hypothesis \ref{Hyp_Orient} or
\ref{Hyp_NonorientFin} or \ref{Hyp_NonorientInf} and satisfies
Hypothesis \ref{Hyp_Assembly5}.  Then there is an injection
\[
\nncc^s_+: \cS_{\TOP +}^B(Y) \longra \cS_\TOP^\Split(Y;Y_0) \x
\UNil_6^s
\]
such that composition with projection onto the first factor
restricts to a subset inclusion $\cS_{\TOP+}^\Split(Y;Y_0) \subseteq
\cS_\TOP^\Split(Y;Y_0)$, and composition with projection onto the
second factor is the manifold splitting obstruction $\Split_L$.
Furthermore, $g$ and $\nncc^s(g)$ have equal image in $\cN_\TOP(Y)$.
\end{enumerate}
\end{thm}

\subsection{Proof by cobordism}

We simply extend Cappell's modification \cite[Chapter
V]{CappellSplit} of the Cappell--Shaneson proof \cite[Theorems 4.1,
5.1]{CS1} of 5-dimensional splitting as to include the non-vanishing
of $\UNil_6^s$. Our homological conditions eschew the performance of
surgery on the 4-manifold $Y_0$. Examples are given in \S
\ref{Sec_SplitFiber}.

\begin{rem}\label{Rem_K}
Friedhelm Waldhausen had shown that $\wt{\Nil}_0$ is a summand of
$\Wh_1(\Pi)$ and that there is an exact sequence of abelian groups
\cite{WaldhausenKrings}:
\begin{multline*}
\Wh_1(\Pi_0) \xra{i_- - i_+} \Wh_1(\Pi_-) \oplus \Wh_1(\Pi_+)
\xra{j_- + j_+}
\Wh_1(\Pi)/\wt{\Nil}_0\\
\xra{\en\bdry\en} \Wh_0(\Pi_0) \xra{i_- - i_+} \Wh_0(\Pi_-) \oplus
\Wh_0(\Pi_+) \xra{j_- + j_+} \Wh_0(\Pi)
\end{multline*}
or, respectively,
\begin{multline*}
\Wh_1(\Pi_0) \xra{i_- - i_+} \Wh_1(\Pi_\infty) \xra{j_\infty}
\Wh_1(\Pi)/\wt{\Nil}_0\\
\xra{\en\bdry\en} \Wh_0(\Pi_0) \xra{i_- - i_+} \Wh_0(\Pi_\infty)
\xra{j_\infty} \Wh_0(\Pi).
\end{multline*}
Waldhausen showed that the cellular splitting obstruction is
algebraically defined as the image $\Split_K(g;Y_0) \in \Wh_0(\Pi_0)
\oplus \wt{\Nil}_0$ of the Whitehead torsion $\tau(g) \in
\Wh_1(\Pi)$. It vanishes if and only if $g$ is CW splittable along
$Y_0$ \cite[erratum]{WaldhausenNotes}.
\end{rem}

\begin{rem}\label{Rem_L}
Sylvain Cappell had shown that $\UNil_6^s$ is a summand of
$L_6^B(\Pi)$ and that there is an exact sequence of abelian groups
\cite{CappellUnitary}:
\begin{multline*}
L_6^h(\Pi_0) \xra{i_- - i_+} L_6^h(\Pi_-) \oplus L_6^h(\Pi_+)
\xra{j_- + j_+}
L_6^B(\Pi) / \UNil_6^s\\
\xra{\en\bdry\en} L_5^h(\Pi_0) \xra{i_- - i_+} L_5^h(\Pi_-) \oplus
L_5^h(\Pi_+) \xra{j_- + j_+} L_5^B(\Pi)
\end{multline*}
or, respectively,
\begin{multline*}
L_6^h(\Pi_0) \xra{i_- - i_+} L_6^h(\Pi_\infty) \xra{j_\infty}
L_6^B(\Pi) / \UNil_6^s\\
\xra{\en\bdry\en} L_5^h(\Pi_0) \xra{i_- - i_+} L_5^h(\Pi_\infty)
\xra{j_\infty} L_5^B(\Pi).
\end{multline*}
If the cellular splitting obstruction vanishes, then Cappell showed
that the manifold splitting obstruction is algebraically defined as
$\Split_L(g;Y_0) \in \UNil_6^s$. It vanishes if $g$ is $\CAT$
splittable along $Y_0$ \cite{CappellFree}.  We shall investigate the
converse.
\end{rem}

\begin{proof}[Proof of Theorem \ref{Thm_Split5}]
(Necessity) Suppose $g$ is $\CAT$ splittable along $Y_0$.  Then
$\Split_K(g;Y_0)=0$ and $\Split_L(g;Y_0)=0$ vanish by Remarks
\ref{Rem_K} and \ref{Rem_L}.

(Sufficiency) Suppose $\Split_K(g;Y_0) = 0$ and $\Split_L(g;Y_0)=0$.
Then $g$ is CW splittable along $Y_0$ and $g \in \cS_\DIFF^B(Y)$
(resp.~$g \in \cS_{\TOP+}^B(Y)$) by Remark \ref{Rem_K}. Since $Y_0$
satisfies the hypotheses in Sections
\ref{Sec_Assembly5}---\ref{Sec_Exactness} for exactness of the
$\CAT$ surgery sequence, by Theorem \ref{Thm_nncc}, it follows that
$\nncc^s(g) = (g,0)$. In other words, $g$ is $\CAT$ splittable along
$Y_0$.

Furthermore, the normal bordisms over $Y_0$ in the proof of Theorem
\ref{Thm_nncc} depend only on the homotopy self-equivalences and
normal self-bordisms of \cite[Proposition 3.5]{Khan_Smoothable4} and
Section \ref{Sec_Assembly5}. Therefore $g: W\to Y$ is $\CAT$
normally bordant to a split homotopy equivalence $g' = g^4$ such
that the $\CAT$ restriction $g': (g')\inv(Y_0) \to Y_0$ is a
homotopy self-equivalence.
\end{proof}

\begin{proof}[Proof of Theorem \ref{Thm_nncc}]
(Definition, I) Let $g: (W,\bdry W) \to (Y,\bdry Y)$ be a homotopy
equivalence with Whitehead torsion $\tau(g) \in B$ and $\bdry g$ a
$\CAT$ isomorphism, representing an element of $\cS_\CAT^B(Y)$. Our
principal goal is to define a $\CAT$ normal bordism $G'$ over $Y \x
\Delta^1$ from $g$ to a homotopy equivalence $g': (W',\bdry W') \to
(Y,\bdry Y)$ such that $h$ is $\CAT$ split along $Y_0$ and that $G'$
has surgery obstruction
\[
\sigma_*(G') \in \UNil_6^s \subseteq L_6^B(\Pi).
\]
Define
\[
\nncc^s(g) := (g',\sigma_*(G')) \in \cS_\CAT^\Split(Y;Y_0) \x
\UNil_6^s.
\]

(Well-definition; Projection properties) Note that $\sigma_*(G')$
depends only on the normal bordism class of $G'$ relative to $\bdry
G' = g \sqcup g'$, and that $\sigma_*(G')$ lies in $L_6^B$ since
$\tau(g), \tau(g') \in B$. Let $Z := \CP^4 \# 2(S^3 \x S^5)$ be the
closed $\CAT$ 8-manifold with Euler characteristic $\chi(Z) = 1$ and
signature $\sigma^*(Z)=1$ used by Weinberger for decorated
periodicity \cite{WeinbergerFibering}. Cappell has shown
\[
\sigma_*(G' \x \Id_Z) = \Split_L(g \x \Id_Z; Y_0 \x Z)
\]
for 13-dimensional homotopy equivalences \cite{CappellSplit}. Note
$\sigma_*(G') = \sigma_*(G' \x \Id_Z)$, by Kwun--Szczarba's torsion
product formula and Ranicki's surgery product formula
\cite[Prop.~8.1(ii)]{RanickiAlgebraicII}.  Also note $\Split_L(g \x
\Id_Z;Y_0 \x Z) = \Split_L(g;Y_0)$, since these splitting
obstructions in $\UNil_6^s$ coincide
\cite[Prop.~7.6.2A]{RanickiExact} with the codimension-one quadratic
signatures \cite[Prop.~7.2.2]{RanickiExact} of $g$ and $g \x \Id_Z$
in the codimension-one Poincar\'e embedding groups $LS_4$ and
$LS_{12}$, and since $\x \Id_Z: LS_4 \to LS_{12}$ is an isomorphism
\cite[Cor.~11.6.1]{Wall}. So $\sigma_*(G') = \Split_L(g;Y_0)$.
Suppose $G''$ is another such $\CAT$ normal bordism from $g$ to some
split $g''$. Then $\sigma_*(G'') = \Split_L(g;Y_0)$. So $G' \cup_g
-G''$ is a normal bordism from $g'$ to $g''$ with surgery
obstruction $0 \in L_6^B$. By the 6-dimensional $\CAT$ $s$-cobordism
theorem, it follows that $g'$ and $g''$ are $\CAT$ isomorphic.
Therefore $\nncc^s(g) = (g', \sigma_*(G'))$ is well-defined and
satisfies the asserted projection properties.

(Bijectivity) Consider Wall's action \cite[Thm.~5.8]{Wall} of the
abelian group $L_6^B$ on the 5-dimensional structure set
$\cS_\CAT^B(Y)$.  It follows from the properties defining $\nncc^s$
that the restriction of Wall's action is the inverse function of
$\nncc^s$:
\[
\act: \cS_\CAT^\Split(Y;Y_0) \x \UNil_6^s \longra \cS_\CAT^B(Y).
\]

(Definition, II) \emph{First, we normally\footnote{The stable normal
$\CAT$ microbundle on $Y$ and on $Y_0$ is $(\ol{g})^*(\nu_W)$ and
$(\ol{g})^*(\nu_W)|Y_0$, lifting the Spivak normal spherical
fibration on $Y$, where $\ol{g}: Y \to W$ is a homotopy inverse of
$g: W \to Y$.} cobord $g$ to a degree one normal map $g^1$ so that
the restriction to $(g^1)\inv(Y_0)$ is a homotopy equivalence.} By
general position, we may assume that $g: W \to Y$ is $\DIFF$
transversal to $Y_0$. Consider the degree one normal map
\[
g_0 := g|_{W_0} : W_0 \longra Y_0
\]
where the $\DIFF$ 4-manifold $W_0 := g\inv(Y_0)$ is the transverse
inverse image of $Y_0$. Denote $\wh{Y} = \wh{Y}_- \cup_{Y_0}
\wh{Y}_+$ as the based cover of $Y$ corresponding to the subgroup
$\Pi_0$.  The $\Z[\Pi_0]$-submodule $P := K_6(\wh{W}_- \x Z)$ is a
finitely generated, projective lagrangian of the
$\Z[\Pi_0]$-equivariant intersection form on the surgery kernel
\[
K_6(W_0 \x Z) = \Ker\prn{g \x \Id_Z: H_6(W_0 \x Z) \longra H_6(Y_0
\x Z)},
\]
where we can homotope $g \x \Id_Z: W \x Z \to Y \x Z$ so that
$g|_{g\inv(Y_0 \x Z)}$ is 6-connected, degree one, normal map
between 12-dimensional manifolds \cite[Lemma II.1]{CappellSplit}.
Furthermore, the projective class $[P] \in \Wh_0(\Pi_0)$ satisfies
$[P^*] = -[P]$ \cite[Lemma II.2]{CappellSplit} and is the
homomorphic image of the Whitehead torsion $\tau(g) = \tau(g \x
\Id_Z) \in \Wh_1(\Pi)$ under Waldhausen's connecting map $\bdry$
(see Remark \ref{Rem_K}).  But $\tau(g) \in B$ implies that $[P] =
\bdry(\tau(g)) = 0$. Therefore the stably free surgery obstruction
vanishes by decorated periodicity:
\[
\sigma_*(g_0) = \sigma_*(g_0 \x \Id_Z) = 0 \in L_4^h(\Pi_0).
\]
Then, by Theorem \ref{Thm_SES}, there exists a $\CAT$ normal bordism
$G^1_0$ from $g_0$ to a homotopy equivalence $g^1_0: W^1_0 \to Y_0$.
So the union
\[
G^1 := (g \x [0,1]) \cup_{g_0 \x 1 \x [-1,1]} (G^1_0 \x [-1,1])
\]
is a $\CAT$ normal bordism from the homotopy equivalence $g: W \to
Y$ to a degree one normal map $g^1: W^1 \to Y$ with transversal
restriction $g^1_0 = (g^1)|_{(g^1)\inv(Y_0)}$ a homotopy
equivalence.

\emph{Second, we normally cobord $g^1$ relative $(g^1)\inv (Y-Y_0)$
to a degree one normal map $g^2$ so that the restriction to
$(g^2)\inv (Y-Y_0)$ has vanishing surgery obstruction.} Since
$g^1_\pm = (g^1)|_{(g^1)\inv(Y_- \sqcup Y_+)}$ (resp.~$g^1_\infty =
(g^1)|_{(g^1)\inv(Y_\infty)}$) restricts to a homotopy equivalence
$g^1_0 \x \set{-1,1}$ on the boundary, and since $G^1$ is a normal
bordism of source and target from $g^1_\pm$ (resp.~$g^1_\infty$) to
the $B$-torsion homotopy equivalence $g$ over the reference space
$K(\Pi,1)$, the image of surgery obstruction vanishes \cite[\S
9]{Wall}:
\[
(j_- + j_+) (\sigma_*(g^1_\pm)) = 0 \in L_5^B(\Pi)
\]
or, respectively,
\[
(j_\infty) (\sigma_*(g^1_\infty)) = 0 \in L_5^B(\Pi).
\]
Therefore, by Cappell's $L$-theory Mayer--Vietoris exact sequence
(Remark \ref{Rem_L}), there exists $a \in L_5^h(\Pi_0)$ such that
\[
(i_- - i_+)(a) = -\sigma_*(g^1_\pm) \in L_5^h(\Pi_-) \oplus
L_5^h(\Pi_+)
\]
or, respectively,
\[
(i_- - i_+)(a) = -\sigma_*(g^1_\infty) \in L_5^h(\Pi_\infty).
\]
Then, by Theorem \ref{Thm_Assembly5} if $\CAT=\TOP$, or by
Hypothesis \ref{Hyp_Surj5} if $\CAT=\DIFF$, there exists a $\CAT$
normal bordism $G^2_0$ from the homotopy equivalence $g^1_0$ to
itself realizing this surgery obstruction: $\sigma_*(G^2_0) = a \in
L_5^h(\Pi_0)$. So the union
\[
G^2 := (g^1 \x [0,1]) \cup_{g^1_0 \x 1 \x [-1,1]} (G^2_0 \x [-1,1])
\]
is a $\CAT$ normal bordism from the degree one normal map $g^1: W^1
\to Y$ to another degree one normal map $g^2: W^2 \to Y$ such that
the transversal restriction $g^2_0 = g^1_0 : W^1_0 \to Y_0$ is a
homotopy equivalence and the transversal restriction $g^2_\pm:
W^2_\pm \to Y_\pm$ (resp.~$g^2_\infty: W^2_\infty \to Y_\infty$) has
vanishing surgery obstruction.

\emph{Third, we normally cobord $g^2$ relative $(g^2)\inv(Y_0)$ to a
degree one normal map $g^3$ so that the restriction to $(g^3)\inv (Y
- Y_0)$ is also a homotopy equivalence.} Since $\sigma_*(g^2_\pm) =
0 \in L_5^h(\Pi_\pm)$ (resp.~$\sigma_*(g^2_\infty) =0 \in
L_5^h(\Pi_\infty)$), by exactness of the 5-dimensional surgery exact
sequence at the $\CAT$ normal invariants \cite[Thm.~10.3]{Wall},
there exists a $\CAT$ normal bordism $G^3_\pm$ (resp.~$G^3_\infty$)
relative $g_2 \x \set{-1,1}$ from the degree one, $\CAT$ normal map
$g^2_\pm$ (resp.~$g^2_\infty$) to a homotopy equivalence $g^3_\pm$
(resp.~$g^3_\infty$). So the union
\[
G^3 := (g^2 \x [0,1] \x [-1,1]) \cup_{g^2_0 \x [0,1] \x \set{-1,1}}
G^3_\pm
\]
or, respectively,
\[
G^3 := (g^2 \x [0,1] \x [-1,1]) \cup_{g^2_0 \x [0,1] \x \set{-1,1}}
G^3_\infty
\]
is a $\CAT$ normal bordism from the degree one normal map $g^2: W^2
\to Y$ to homotopy equivalence $g^3 = g^3_- \cup_{g^3_0} g^3_+: W^3
\to Y$ (resp.~$g^3 = \cup_{g^3_0}\, g^3_\infty: W^3 \to Y$) which is
$\CAT$ split along $Y_0$.

\emph{Finally, we normally cobord the split homotopy equivalence
$g^3$ to another split homotopy equivalence $g^4 = g'$ so that the
normal bordism $G' = G^1 \cup G^2 \cup G^3 \cup G^4$ from $g$ to
$g'$ has surgery obstruction in the subgroup $\UNil_6^s$ of the
abelian group $L_6^B$.} Consider the surgery obstruction of the
$\CAT$ normal bordism from $g$ to $g^3$:
\[
b := -\sigma_*(G^1 \cup G^2 \cup G^3) \in L_6^B(\Pi).
\]
Let $c := \bdry(b) \in L_5^h(\Pi_0)$ be the image in Cappell's
$L$-theory Mayer--Vietoris exact sequence (Remark \ref{Rem_L}).  By
Theorem \ref{Thm_Assembly5} if $\CAT=\TOP$, or by Hypothesis
\ref{Hyp_Surj5} if $\CAT=\DIFF$, there exists a $\CAT$ normal
bordism $G^{3.5}_0$ from the homotopy equivalence $g^3_0$ to itself
realizing this surgery obstruction: $\sigma_*(G^{3.5}_0) = c \in
L_5^h(\Pi_0)$.  Then
\[
0 = (i_- - i_+)(c) = \sigma_*(g^3_\pm \cup G^{3.5}_0 \x \set{-1,1})
\in L_5^h(\Pi_-) \oplus L_5^h(\Pi_+)
\]
or, respectively,
\[
0 = (i_- - i_+)(c) = \sigma_*(g^3_\infty \cup G^{3.5}_0 \x
\set{-1,1}) \in L_5^h(\Pi_\infty).
\]
Therefore, by exactness of the 5-dimensional surgery exact sequence
at the $\CAT$ normal invariants \cite[Thm.~10.3]{Wall}, there exists
a $\CAT$ normal bordism $G^{3.5}_\pm$ (resp.~$G^{3.5}_\infty$)
relative $g^3_0 \x \set{-1,1}$ from the degree one normal map
$g^3_\pm \cup G^{3.5}_0 \x \set{-1,1}$ (resp.~$g^3_\infty \cup
G^{3.5}_0 \x \set{-1,1}$) to a homotopy equivalence $g^{3.5}_\pm:
W^{3.5}_\pm \to Y_\pm$ (resp.~$g^{3.5}_\infty: W^{3.5}_\infty \to
Y_\infty$). So the union
\[
G^{3.5} := (G^{3.5}_0 \x [-1,1]) \cup_{G^{3.5}_0 \x \set{-1,1}}
G^{3.5}_\pm
\]
or, respectively,
\[
G^{3.5} := (G^{3.5}_0 \x [-1,1]) \cup_{G^{3.5}_0 \x \set{-1,1}}
G^{3.5}_\infty
\]
is a $\CAT$ normal bordism from the split homotopy equivalence $g^3:
W^3 \to Y$ to another split homotopy equivalence $g^{3.5}: W^{3.5}
\to Y$ such that
\[
\sigma_*(G^1 \cup G^2 \cup G^3 \cup G^{3.5}) = j(d) \oplus e \in
L_6^B(\Pi)
\]
for some $d \in L_6^h(\Pi_-) \oplus L_6^h(\Pi_+)$ (resp.~$d \in
L_6^h(\Pi_\infty)$) and $e \in \UNil_6^s$.

By Wall realization on 5-dimensional $\CAT$ structure sets
\cite[Thm.~10.5]{Wall}, there exists a $\CAT$ normal bordism
$G^4_\pm$ (resp.~$G^4_\infty$) relative $g^{3.5}_0 \x \set{-1,1}$
from the homotopy equivalence $g^{3.5}_\pm$ (resp.~$g^{3.5}_\infty$)
to another homotopy equivalence $g^4_\pm$ (resp.~$g^4_\infty$) such
that $\sigma_*(G^4_\pm) = -d$ (resp.~$\sigma_*(G^4_\infty) = -d$).
So the union
\[
G^4 := G^{3.5} \cup_{g^{3.5} \x 0} (g^{3.5}_0 \x [0,1] \x [-1,1]
\cup_{} G^4_\pm)
\]
or, respectively,
\[
G^4 := G^{3.5} \cup_{g^{3.5} \x 0} (g^{3.5}_0 \x [0,1] \x [-1,1]
\cup_{} G^4_\infty)
\]
is a $\CAT$ normal bordism from the split homotopy equivalence
$g^{3.5}: W^{3.5} \to Y$ to another split homotopy equivalence $g^4
: W^4 \to Y$ such that
\[
\sigma_*(G^1 \cup G^2 \cup G^3 \cup G^4) = e \in \UNil_6^s.
\]
Thus the definition of $\nncc^s$ is complete.
\end{proof}

\section{Fibering and splitting over the circle}\label{Sec_SplitFiber}

We approach the problem of fibering a 5-manifold $W$ over the circle
$S^1$ from Farrell's point of view \cite{Farrell_ICM1970}, which
involves a finite domination of the infinite cyclic cover $\ol{W}$,
the covering translation $t: \ol{W} \to \ol{W}$, and a certain
mapping torus.

Gluing the ends of a self $h$-cobordism $(Y;X,X)$ by a self
homeomorphism $\alpha: X \to X$ yields an \textbf{$h$-block bundle
$\cup_\alpha Y$ over $S^1$}
\cite[p.~306]{Chapman_ApproxHilbertCube}. This is classically known
as a \emph{pseudo-fibering over $S^1$
\cite[Defn.~3.1]{Farrell_fibering}
\cite[Defn.~4.2]{SiebenmannFibering}}. Consider the zero-torsion
version of $h$-block bundles over $S^1$.

\begin{defn}\label{Defn_sBlock}
We call $E$ the total space of a \textbf{$\CAT$ $s$-block bundle
over $S^1$ with homotopy fiber $X$} if $E$ is the compact $\CAT$
manifold obtained by gluing the ends of a self $\CAT$ $s$-cobordism
$(Y;X,X) \rel \bdry X$ by a $\CAT$ automorphism $\alpha: X \to X$.
We write $E = \cup_\alpha Y := Y / (x,0) \sim (\alpha(x),1)$, and a
special case is the mapping torus $X \rtimes_\alpha S^1 :=
\cup_\alpha X \x \Delta^1$. The induced continuous map $E \to S^1$
is called a \textbf{$\CAT$ $s$-block bundle projection}, which is
unique up to homotopy.
\end{defn}

Fiber bundles are special cases of $s$-block bundles, and the
converse is true if the $s$-cobordism theorem holds for the fiber in
the given manifold category.

\begin{defn}[{Quinn, compare \cite[\S2.3]{Nicas_Quinn}}]
Let $X$ be a compact $\CAT$ manifold. The \textbf{block homotopy
automorphism space} $\wt{\G}^s(X)$ is the geometric realization of
the Kan $\Delta$-set whose $k$-simplices are simple homotopy
self-equivalences $e: X \x \Delta^k \to X \x \Delta^k$ of
$(k+2)$-ads which restrict to the identity over $\bdry X$. Note that
$\hAut^s(X) = \pi_0 \wt{\G}^s(X)$. The basepoint is the identity
$\Id_X: X \to X$.

Similarly, the \textbf{block structure space} $\wt{\bbS}^s_\CAT(X)$
is the geometric realization of the Kan $\Delta$-set whose
$k$-simplices are simple homotopy equivalences $Y \to X \x \Delta^k$
of $\CAT$ manifold $(k+2)$-ads in $\R^\infty$ which restrict to
$\CAT$ isomorphisms over $\bdry X$. Note the $\CAT$ $s$-bordism
structure set is $\cS^s_\CAT(X) = \pi_0 \wt{\bbS}^s_\CAT(X)$.  We
define the decoration $\CAT=\TOP+$ for block structures to mean that
the $\TOP$ manifolds $Y$ are smoothable, that is, without a
preference of $\DIFF$ structure.

Naturally, there is a simplicial inclusion $\wt{\G}^s(X) \hookra
\wt{\bbS}^s_\CAT(X)$.
\end{defn}

An assembly-type function over $S^1$ is described as follows.

\begin{defn}\label{Defn_cupBlock}
Define an \textbf{$\alpha$-twisted simplicial loop in
$(\wt{\bbS}^s_\CAT(X), \wt{\G}^s(X))$} as a simple homotopy
equivalence $(h;h_0,h_1): (Y;X,X) \to X \x (\Delta^1;0,1)$ of $\CAT$
manifold triads such that the simple homotopy self-equivalences
$h_i: X \to X$ satisfy $h_1 = \alpha \circ h_0$ and that $h$
restricts to a $\CAT$ isomorphism over $\bdry X$.  We define the
\textbf{$\alpha$-twisted fundamental set}
$\pi_1^\alpha(\wt{\bbS}^s_\CAT(X), \wt{\G}^s(X))$ as the set of
homotopy classes of these $\alpha$-twisted loops. Note, if $\alpha$
is the identity automorphism, then this set is the first homotopy
set of the pair. Define the \textbf{union function}
\[
\cup: \pi_1^\alpha(\wt{\bbS}^s_\CAT(X), \wt{\G}^s(X)) \longra \cS_\CAT^s(X \rtimes_\alpha S^1)
\]
as follows.  Let $(h;h_0,h_1): (Y;X,X) \to X \x (\Delta^1;0,1)$ be
an $\alpha$-twisted simplicial loop. Then there is an induced simple
homotopy equivalence, well-defined on homotopy classes of loops:
\[
\cup (h;h_0,h_1): \cup_{\Id_X} Y \longra X \rtimes_\alpha
S^1; \quad [y] \longmapsto [h(y)].
\]
\end{defn}

Now we are ready to study three families of examples of homotopy
fibers.

\subsection{Statement of results}

Our theorems below are crafted as to eliminate any algebraic $K$- or
$L$-theory obstructions to splitting and fibering over the circle.

The first splitting theorem (\ref{Thm_SplittingOverCircle}) is a
special case of the general splitting theorem (\ref{Thm_Split5}).
Here, for any homotopy self-equivalence $h: (X,\bdry X) \to (X,\bdry
X)$, the mapping torus $Y = X \rtimes_h S^1$ is a Poincar\'e pair
with $X$ a two-sided Poincar\'e subpair
\cite[Prop.~24.4]{RanickiKnot}. This level of abstraction is
required to prove the fibering theorem
(\ref{Thm_FiberingOverCircle}).

\begin{thm}\label{Thm_SplittingOverCircle}
Let $X$ be any of the 4-manifolds $Q,S,F$ defined in
(\ref{Hyp_Projective}, \ref{Hyp_Surface}, \ref{Hyp_Sum}). Let $h: X
\to X$ be a homotopy equivalence which restricts to a diffeomorphism
$\bdry h: \bdry X \to \bdry X$. Suppose $M$ is a compact $\DIFF$
5-manifold and $g: M \to X \rtimes_h S^1$ is a homotopy equivalence
which restricts to a diffeomorphism $\bdry g: \bdry M \to \bdry X
\rtimes_{\bdry h} S^1$.

Then $g$ is homotopic to a map $g'$ which restricts to a simple
homotopy equivalence $g': X' \to X$ such that the $\TOP$ inverse
image $X' := (g')\inv(X)$ is homeomorphic to $X$ and the exterior
$M'$ of $X'$ in $M$ is a smoothable $\TOP$ self $s$-cobordism of
$X$.
\end{thm}

The second splitting theorem (\ref{Thm_StructuresOverCircle})
connects homotopy $\TOP$ structures on mapping tori to smoothable
$s$-cobordisms, homotopy self-equivalences, and the stable smoothing
invariant of Kirby and Siebenmann \cite{KS}.

\begin{thm}\label{Thm_StructuresOverCircle}
Let $X$ be any of the 4-manifolds $Q,S,F$ defined in
(\ref{Hyp_Projective}, \ref{Hyp_Surface}, \ref{Hyp_Sum}). Let
$\alpha: X \to X$ be a diffeomorphism.  Then there is an exact
sequence of based sets:
\[
\pi_1^\alpha(\wt{\bbS}^s_{\TOP+}(X), \wt{\G}^s(X)) \xra{\en\cup\en}
\cS_\TOP^h(X \rtimes_\alpha S^1) \xra{\en\ks\en} \F_2 \oplus
H_1(X;\F_2)_\alpha.
\]
\end{thm}

Our fibering theorem (\ref{Thm_FiberingOverCircle}) is proven using
a key strategy of Tom Farrell \cite{Farrell_ICM1970}, except we do
not require 4-dimensional Siebenmann ends on the infinite cyclic
cover $\ol{M}$ to exist.  A \textbf{connected manifold band} $(M,f)$
consists of a connected manifold $M$ and a continuous map $f: M \to
S^1$ such that the infinite cyclic cover $\ol{M} := f^*(\R)$ is
connected (i.e.~$f_*: \pi_1(M) \to \pi_1(S^1)$ surjective) and is
dominated by a finite CW complex \cite[Defn.~15.3]{HughesRanicki}.
Observe that the manifold $\ol{M}$ is a strong deformation retract
of the homotopy fiber $\hofiber(f)$.  If $f: M \to S^1$ is homotopic
to a fiber bundle projection, say with fiber $X'$, then
$\hofiber(f)$ is homotopy equivalent to $X'$.

\begin{thm}\label{Thm_FiberingOverCircle}
Let $X$ be any of the 4-manifolds $Q,S,F$ defined in
(\ref{Hyp_Projective}, \ref{Hyp_Surface}, \ref{Hyp_Sum}). Let
$(M,f)$ be a connected $\DIFF$ 5-manifold band such that $\bdry X
\to \bdry M \xra{\bdry f} S^1$ is a $\DIFF$ fiber bundle and the
homotopy equivalence $\bdry X \to \hofiber(\bdry f)$ extends to a
homotopy equivalence $X \to \hofiber(f)$. Then $f: M \to S^1$ is
homotopic $\!\!\rel \bdry M$ to the projection of a smoothable
$\TOP$ $s$-block bundle with fiber $X$.
\end{thm}

\begin{rem}
Topological splitting and fibering of 5-manifolds $W$ over the
circle $S^1$ with fibers like $T^4$ or $Kl \x S^2$ can be
established by Weinberger's splitting theorem
\cite{WeinbergerFibering}, since the fundamental groups $\Z^4$ and
$\Z \rtimes \Z$ have subexponential growth.  A fortiori,
simply-connected topological 4-manifolds $X$ are allowable fibers in
Weinberger's fibering theorem. These 4-manifolds $X$ are classified,
by Milnor in the $\PL$ case and Freedman--Quinn in the $\TOP$ case
\cite{FQ}, up to homotopy equivalence by their intersection form.
Unfortunately, the smooth splitting and fibering problems for such
$W$ remain unsolved, even as $\DIFF$ $s$-block bundle maps.
\end{rem}

Now let us state the smooth result promised in the Introduction.

\begin{thm}\label{Thm_SmoothFiberingOverCircle}
Consider the closed, non-orientable, smooth 4-manifolds $Q$
(\ref{Hyp_Projective}). Let $(M,f)$ be a connected $\DIFF$
5-manifold band such that $Q$ is homotopy equivalent to
$\hofiber(f)$. Then $f: M \to S^1$ is homotopic to the projection of
a $\DIFF$ $s$-block bundle with fiber $Q$.
\end{thm}

\subsection{Vanishing of lower Whitehead groups}

We start by showing that every homotopy equivalence under
consideration has zero Whitehead torsion.  Comparable results are
found by J.~Hillman \cite[\S 6.1]{Hillman4Book} and intensely use
\cite[\S 19]{WaldhausenKrings}.

\begin{prop}\label{Prop_VanishingWhiteheadGroup}
Let $X$ be any of the  4-manifolds $Q,S,F$ defined in
(\ref{Hyp_Projective}, \ref{Hyp_Surface}, \ref{Hyp_Sum}). Let $h: X
\to X$ be a homotopy self-equivalence which restricts to a
diffeomorphism $\bdry h: \bdry X \to \bdry X$. Suppose $Y = X
\rtimes_h S^1$ is the total space of a mapping torus fibration $X
\to Y \to S^1$ with monodromy $h$. Then the Whitehead groups vanish
for all $* \leq 1$:
\[
\Wh_*(\pi_1 X) = \Wh_*(\pi_1 Y) = 0.
\]
\end{prop}

Note if $h$ is a diffeomorphism, then $Y$ is the total space of a
$\DIFF$ fiber bundle.

\begin{proof}
Each case of fiber shall be handled separately.

\indent\emph{Case of fiber $Q$ (\ref{Hyp_Projective}). } Let $Q$ be
any of the non-orientable 4-manifolds listed. There are no
non-identity automorphisms of the cyclic group $\pi_1(Q) = C_2$ of
order two, so $\pi_1(Y) = C_2 \x C_\infty$.  S.~Kwasik has shown
that $\Wh_1(\pi_1 Q) = \Wh_1(\pi_1 Y) = 0$, using the Rim square for
the ring $\Z[C_\infty][C_2]$, the Bass--Heller--Swan decomposition,
and assorted facts \cite[pp.~422--423]{KwasikLowDimScob}.

\indent\emph{Case of fiber $S$ (\ref{Hyp_Surface}). } There are
three types of $S$.

\emph{Type 1.} Suppose $S$ is the total space of a fiber bundle $S^2
\to S \to \Sigma$ such that $\Sigma$ is a compact, connected,
possibly non-orientable 2-manifold of positive genus. Then, by
\cite[Theorem 19.5(5)]{WaldhausenKrings}, the fundamental group
$\pi_1(S) = \pi_1(\Sigma)$ is a member of Waldhausen's class $\Cl$
of torsion-free groups.  So, by \cite[Proposition
19.3]{WaldhausenKrings}, the HNN-extension $\pi_1(Y) = \pi_1(S)
\rtimes C_\infty$ is a member of $\Cl$. Therefore, by \cite[Theorem
19.4]{WaldhausenKrings}, we obtain that $\Wh_*(\pi_1 S) =
\Wh_*(\pi_1 Y)=0$.

\emph{Type 2.} Suppose $S$ is the total space of a fiber bundle $H
\to S \to S^1$ such that $H$ is a closed, connected, hyperbolic
3-manifold.  By Mostow rigidity, we may select the monodromy
diffeomorphism $H \to H$ to be an isometry \cite{MostowSmooth} up to
smooth isotopy \cite{GMT}. Since $H \to S \to S^1$ has isometric
monodromy implies that $S$ is isometrically covered by $\qH^3 \x
\qE^1$, the curvature matrix is constant, hence $S$ is $A$-regular.
By hypothesis, $h$ is homotopic rel $\bdry S$ to an isometry of
$S-\bdry S$. Then $Y-\bdry Y$ is isometrically covered by $\qH^3 \x
\qE^2$, hence $Y-\bdry Y$ is $A$-regular. Therefore, by \cite[Lemma
0.12]{FJ_TorsionFreeGL}, we obtain that $\Wh_*(\pi_1 S) =
\Wh_*(\pi_1 Y) = 0$ for all $* \leq 1$.

\emph{Type 3.} Suppose the interior $S-\bdry S$ admits a complete,
finite volume metric of euclidean type (resp.~real hyperbolic or
complex hyperbolic). Since $S-\bdry S$ is isometrically covered by
$\qE^4$ (resp.~$\qH^4$ or $\CH^2$), the curvature matrix is
constant, hence $S-\bdry S$ is $A$-regular. By hypothesis, $h$ is
homotopic rel $\bdry S$ to an isometry of $S-\bdry S$. Then $Y -
\bdry Y$ is isometrically covered by $\qE^4 \x \qE^1$ (resp.~$\qH^4
\x \qE^1$ or $\CH^2 \x \qE^1$), the curvature matrix is constant,
hence $Y - \bdry Y$ is $A$-regular. Therefore, by \cite[Lemma
0.12]{FJ_TorsionFreeGL}, we obtain that $\Wh_*(\pi_1 S) =
\Wh_*(\pi_1 Y) = 0$ for all $* \leq 1$.

\indent\emph{Case of fiber $F$ (\ref{Hyp_Sum}). } There are two
types of connected summands $F_i$.

\emph{Type 1. } Suppose $F_i$ is the total space of a fiber bundle
$H_i \to F_i \to S^1$ such that the compact, connected, irreducible,
orientable 3-manifold $H_i$ either is $S^3$ or $D^3$ or has non-zero
first Betti number. Then, by \cite[Proposition
19.5(6,8)]{WaldhausenKrings}, we obtain that $\pi_1(H_i)$ is a
member of $\Cl$. So, by \cite[Proposition 19.3]{WaldhausenKrings},
the HNN-extension $\pi_1(F_i) = \pi_1(H_i) \rtimes C_\infty$ is a
member of $\Cl$.

\emph{Type 2. } Suppose $F_i$ is the total space of a fiber bundle
$\Sigma_i^f \to F_i \to \Sigma_i^b$ such that the fiber and base are
compact, connected, orientable 2-manifolds of positive genus. Then,
by a theorem of J.~Hillman \cite[Thm.~1]{Hillman_SurfaceBundles}, we
obtain that $\pi_1(F_i)$ is a member of $\Cl$. Indeed, the direct
algebraic proof of Cavicchioli--Hegenbarth--Spaggiari
\cite[Thm.~3.12]{CHS} uses a Mayer--Vietoris argument for a
connected-sum decomposition of the base $\Sigma_i^b$, which extends
to aspherical, compact, possibly non-orientable surfaces with
possibly non-empty boundary.

\emph{Conclusion. } Now, since $\pi_1(F_i)$ is a member of $\Cl$, by
\cite[Proposition 19.3]{WaldhausenKrings}, the fundamental group
$\pi_1(F) = \bigstar_{i=1}^n \pi_1(F_i)$ of the connected sum $F =
F_1 \# \cdots \# F_n$ is a member of $\Cl$.  So the HNN-extension
$\pi_1(Y) = \pi_1(F) \rtimes C_\infty$ is a member of $\Cl$.
Therefore, by \cite[Theorem 19.4]{WaldhausenKrings}, we obtain
$\Wh_*(\pi_1 F) = \Wh_*(\pi_1 Y)=0$.
\end{proof}

\subsection{Proof of main theorems over the circle}

\begin{proof}[Proof of Theorem \ref{Thm_SplittingOverCircle}]
By the general splitting theorem (Thm.~\ref{Thm_Split5}), it
suffices to show that the following conditions hold, in order:
\begin{enumerate}
\item
$X$ satisfies Hypothesis \ref{Hyp_Orient} or
\ref{Hyp_NonorientFin} or \ref{Hyp_NonorientInf} and Hypothesis
\ref{Hyp_Assembly5},

\item
the obstructions $\Split_K(g;X)$ and $\Split_L(g;X)$ vanish, and

\item
the $h$-cobordism $(M';X',X')$ and homotopy equivalence $g': X' \to
X$ have zero Whitehead torsion.
\end{enumerate}

\emph{Condition (1).} \emph{Case of $X=Q$.} By Corollary
\ref{Cor_OrderTwo}, $Q$ satisfies Hypothesis \ref{Hyp_NonorientFin}.
Since \cite[Theorem 13A.1]{Wall} implies $L_5^h(\Z[C_2]^-)=0$, $Q$
fulfills Theorem \ref{Thm_Assembly5}.

\emph{Case of $X=S$.} There are three types of $S$.  If $S$ is
non-orientable, then, since the abelianization $H_1(S;\Z)$ is
2-torsionfree, there exists a lift $\hat{\omega}: \pi_1(S) \to \Z$
of the orientation character $\omega: \pi_1(S) \to \Z^\x$.

\emph{Type 1.} By Corollary \ref{Cor_Assembly5Nonaspherical}, $S$
satisfies Hypothesis \ref{Hyp_Assembly5}.  If $S$ is orientable,
then, by Corollary \ref{Cor_FreeProdCrys}, $S$ satisfies Hypothesis
\ref{Hyp_Orient}. Otherwise, if $S$ is non-orientable, then, by
Corollary \ref{Cor_Klein}, $S$ satisfies Hypothesis
\ref{Hyp_NonorientInf}.

\emph{Type 2.} Since $H$ is irreducible, by Corollary
\ref{Cor_Assembly5}, $S$ satisfies Hypothesis \ref{Hyp_Assembly5}.
Recall that $S$ can be isometrically covered by $\qH^3 \x \qE^1$, by
Mostow rigidity \cite{MostowSmooth, GMT}. Then, by \cite[Proposition
0.10]{FJ_TorsionFreeGL} and the argument of Corollary
\ref{Cor_FreeProdCrys}, it follows that $\kappa_2$ is injective. If
$S$ is orientable, then $S$ satisfies Hypothesis \ref{Hyp_Orient}.
Otherwise, if $S$ is non-orientable, then, since the lift
$\hat{\omega}$ exists, $S$ satisfies Hypothesis
\ref{Hyp_NonorientInf}.

\emph{Type 3.} By Corollary \ref{Cor_Assembly5}, $S$ satisfies
Hypothesis \ref{Hyp_Assembly5}. If $S$ is orientable, then, by
Corollary \ref{Cor_FreeProdCrys}, $S$ satisfies Hypothesis
\ref{Hyp_Orient}. Otherwise, suppose $S$ is non-orientable. The
argument in Corollary \ref{Cor_FreeProdCrys} for the injectivity of
$\kappa_2$ is the same as if $S$ were orientable. Then, since
$\hat{\omega}$ exists, $S$ satisfies Hypothesis
\ref{Hyp_NonorientInf}.

\emph{Case of $X=F$.} There are two types of connected summands
$F_i$. By Corollary \ref{Cor_Assembly5}, $F_i$ satisfies Hypothesis
\ref{Hyp_Assembly5}. Note that $F_i$ satisfies Hypothesis
\ref{Hyp_Orient}, by Corollary \ref{Cor_Haken} for \emph{Type 1} and
by Corollary \ref{Cor_Klein} for \emph{Type 2}.

\emph{Condition (2).} Recall the notation $\Pi_0 := \pi_1(X)$ and
$\Pi := \pi_1(X \rtimes_h S^1) \iso \pi_1(M)$. The $K$-theory
obstruction $\Split_K(g;X)$ lies in
\[
\Wh_0(\Pi_0) \oplus \wt{\Nil}_0(\Z[\Pi_0];0,0,{}_- \Z[\Pi_0]_+, {}_+
\Z[\Pi_0]_-).
\]
By \cite[Theorem 2]{WaldhausenKrings}, the second factor is a
summand of $\Wh_1(\Pi)$.  Both $\Wh_0(\Pi_0)$ and $\Wh_1(\Pi)$
vanish by Proposition \ref{Prop_VanishingWhiteheadGroup}. The
$L$-theory obstruction $\Split_L(g;X)$ lies in
$\UNil_6^s(\Z[\Pi_0];0,0)$, which vanishes by definition \cite[\S
1]{CappellUnitary}.

\emph{Condition (3).} The torsions of the $h$-cobordism $(M';X',X')$
and the homotopy equivalence $g': X' \to X$ lie in $\Wh_1(\Pi_0)$,
which vanishes by Proposition \ref{Prop_VanishingWhiteheadGroup}.
\end{proof}

\begin{proof}[Proof of Theorem \ref{Thm_StructuresOverCircle}]
Let $[g] \in \cS_\TOP^h(X \rtimes_\alpha S^1)$. Then $[g]$ is an
$h$-bordism class of some homotopy equivalence $g: M \to X
\rtimes_\alpha S^1$ such that $M$ is a compact $\TOP$ 5-manifold and
the restriction $\bdry g: \bdry M \to \bdry X \rtimes_{\bdry \alpha}
S^1$ is a homeomorphism. Since $\alpha: X\to X$ is a diffeomorphism
implies that the mapping torus $X \rtimes_\alpha S^1$ is a $\DIFF$
5-manifold, we have
\[
\ks[g] = g_*\ks(M) \in \F_2 \oplus H_1(X; \F_2)_\alpha \iso H^4(X \rtimes_\alpha
S^1, \bdry X \rtimes_{\bdry \alpha} S^1; \F_2).
\]
Here, the isomorphism is obtained from Poincar\'e duality and the
Wang sequence. Note, by definition, that the basepoints are
respected by $\cup$ and $\ks$, and that the composite $\ks \circ
\cup$ vanishes.

Suppose $\ks[g]=0$. Since $\bdry g$ induces a $\DIFF$ structure on
$\bdry M$ and $\ks(M)=0$ and $\dim(M)>4$, by a consequence
\cite[Thm.~IV.10.1]{KS} of Milnor's Lifting Criterion and the
Product Structure Theorem, the $\DIFF$ structure on $\bdry M$
extends to a $\DIFF$ structure on $M$. So, by Theorem
\ref{Thm_SplittingOverCircle}, we obtain that $g$ is homotopic to a
$\TOP$ split homotopy equivalence $g': M \to X \rtimes_\alpha S^1$
such that the restriction $g': (g')\inv(X) \to X$ is a homotopy
self-equivalence. Moreover, the restriction of $g'$ to the exterior
of $X$ yields a smoothable $\TOP$ $s$-cobordism $(M';X,X)$ and an
$\alpha$-twisted simplicial loop $(g'_\infty;g'_0,g'_1): (M';X,X)
\to X \x (\Delta^1;0,1)$ in $(\wt{\bbS}^s_{\TOP+}(X),
\wt{\G}^s(X))$. Therefore $[g] = [g'] = \cup [g'_\infty,
g'_0,g'_1]$. Thus exactness is proven at $\cS_\TOP^h(X
\rtimes_\alpha S^1)$.
\end{proof}

\begin{proof}[Proof of Theorem \ref{Thm_FiberingOverCircle}]
The hypothesis gives a homotopy equivalence $d: (X,\bdry X) \to
(\ol{M}, \bdry \ol{M})$ of pairs with homotopy inverse $u: (\ol{M},
\bdry \ol{M}) \to (X, \bdry X)$ such that $\bdry u \circ \bdry d =
\Id_{\bdry X}$. In particular, $d$ is a domination of $(\ol{M},
\bdry \ol{M})$ by $(X, \bdry X)$. That is, there is a homotopy $H:
[0,1] \x \ol{M} \to \ol{M}$ such that $H_0 = d \circ u$ and $H_1 =
\Id_{\ol{M}}$.

Recall $\ol{M} = f^*(\R)$. Let $\ol{f}: \ol{M} \to \R$ be the
sub-projection covering $f: M \to S^1$. Let $t: \ol{M} \to \ol{M}$
be the unique covering transformation such that $\ol{f} t(x) =
\ol{f}(x) + 1$. Then the following composite is a homotopy
self-equivalence of pairs:
\[
h := u \circ t \circ d: (X, \bdry X) \longra (X, \bdry X).
\]
So, by cellular approximation of $h$ and \cite[Proposition
24.4]{RanickiKnot}, the mapping torus $X \rtimes_h S^1 = X \x
[0,1]/(x,0)\sim (h(x),1)$ is a finite Poincar\'e pair of formal
dimension 5, such that $X \iso X \x [1/2]$ is a two-sided Poincar\'e
subpair with tubular neighborhood $X \x [-1,1] \iso X \x [1/3,2/3]$.
Furthermore, since $\bdry d: \bdry X \to \bdry \ol{M}$ and $\bdry u:
\bdry \ol{M} \to \bdry X$ are the $0$-section and projection from
$\bdry \ol{M} = \bdry X \x \R$, we obtain that the homotopy
self-equivalence $\bdry h: \bdry X \to \bdry X$ is in fact a
self-diffeomorphism on each connected component. In particular,
$\bdry M = \bdry X \rtimes_{\bdry h} S^1$.

Observe that the Borel construction fits into a fiber bundle $\R \to
\ol{M} \x_{\Z} \R \to M$ and similarly for $\bdry M$. Then the
projection $g_1: \ol{M} \x_{\Z} \R = \ol{M} \rtimes_t S^1 \to M$ is
a homotopy equivalence of manifold pairs. Note $\set{H_s \circ t
\circ d}_{s \in [0,1]}$ is a homotopy from $d \circ h$ to $t \circ
d: X \to \ol{M}$. Define a continuous map
\[
g_2: X \rtimes_h S^1 \longra \ol{M} \rtimes_t S^1; \quad [x,s]
\longmapsto [H(s,td(x)),s].
\]
By cyclic permutation of the composition factors of $h$, and by the
adjunction lemma (see \cite{SiebenmannFibering}), the map $g_2$ is a
homotopy equivalence of manifold pairs.  Let $\ol{g}_i$ be a
homotopy inverse of $g_i$ for $i=1,2$. Then we obtain a homotopy
equivalence
\[
g := \ol{g}_2 \circ \ol{g}_1 : (M, \bdry M) \longra (X \rtimes_h S^1, \bdry X \rtimes_{\bdry h} S^1).
\]
Furthermore, since $\bdry X \to \bdry M \to S^1$ is already a fiber
bundle, the homotopy inverse $\bdry \ol{g} = \bdry g_1 \circ \bdry
g_2$ is homotopic to the above diffeomorphism $\bdry X
\rtimes_{\bdry h} S^1 \to \bdry M$. By Theorem
\ref{Thm_SplittingOverCircle}, the homotopy equivalence $g$ is
homotopic $\!\!\rel \bdry M$ to a map $g'$ such that the $\TOP$
transverse restriction $g': X' := (g')\inv(X) \to X$ is a simple
homotopy equivalence and there is a homeomorphism $X' \homeo X$.
Moreover, $M = \cup_{\Id_{X'}} M'$ is obtained by gluing the ends of
the smoothable $\TOP$ self $s$-cobordism $M' := M - X' \x (-1,1)$ by
the identity map.

Define quotient maps
\begin{eqnarray*}
q: X \rtimes_h S^1 \longra S^1; && [x,s] \longmapsto [s]\\
q': M \longra S^1; && q' := q \circ g'.
\end{eqnarray*}
Note $\bdry q' = \bdry f: \bdry M \to S^1$ is the fiber bundle
projection. Therefore, by obstruction theory, the continuous map $f:
M \to S^1$ and the $\TOP$ $s$-block bundle projection $q': M \to
S^1$ are homotopic $\!\!\rel \bdry M$ if and only if they determine
the same kernel subgroup of $\pi_1(M)$.  Then, by covering space
theory, it suffices to show that the isomorphism $g_*: \pi_1(M) \to
\pi_1(X \rtimes_h S^1)$ maps the subgroup $\Ker(f_*) =
p_*\pi_1(\ol{M})$ onto the subgroup $\Ker(q_*) = p'_*\pi_1(X\x \R)$.
Here, $p: \ol{M} \to M$ and $p': X\x \R \to X \rtimes_h S^1$ are the
infinite cyclic covers. Observe that the $\pi_1$-isomorphism induced
by the split homotopy equivalence $g_2: X \rtimes_h S^1 \to \ol{M}
\rtimes_t S^1$ maps the subgroup $\Ker(q_*) = \pi_1(X)$ onto
$\pi_1(\ol{M})$, and that the $\pi_1$-isomorphism induced by the
homotopy equivalence $g_1: \ol{M} \rtimes_t S^1 \to M$ maps the
subgroup $\pi_1(\ol{M})$ onto $\Ker(f_*)$. So, since $g_1 \circ g_2
= \ol{g}$ is the homotopy inverse of $g$, we are done.
\end{proof}

\begin{proof}[Proof of Theorem \ref{Thm_SmoothFiberingOverCircle}]
The proof of Theorem \ref{Thm_FiberingOverCircle} constructs
homotopy equivalences $h: Q \to Q$ and $g: M \to Q \rtimes_h S^1$.
Observe Corollary \ref{Cor_OrderTwo} implies that $Q$ satisfies
Hypothesis \ref{Hyp_NonorientFin}, and Remark \ref{Rem_OrderTwo}
implies that $Q$ satisfies Hypothesis \ref{Hyp_Surj5}. Recall that
Conditions (2) and (3) of Proof \ref{Thm_SplittingOverCircle} hold.
Then, by Theorem \ref{Thm_Split5}, the homotopy equivalence $g$ is
homotopic to a map $g'$ such that the $\DIFF$ transverse restriction
$g': Q' := (g')\inv(Q) \to Q$ is a simple homotopy equivalence and
there is a diffeomorphism $Q' \homeo Q$. Moreover, the $\DIFF$
5-manifold $M = \cup_{\Id_{Q'}} M'$ is obtained by gluing the ends
of the $\DIFF$ self $s$-cobordism $M' := M - Q' \x (-1,1)$ by the
identity map. The remainder of Proof \ref{Thm_FiberingOverCircle}
shows that $f: M \to S^1$ is homotopic to the $\DIFF$ $s$-block
bundle projection $q': \cup_{\Id_{Q}} (M';Q,Q) \to S^1$ obtained
from $g'$.
\end{proof}
\subsection*{Acknowledgments}

The author is grateful to Jim Davis, who supervised the splitting
result of this paper, which was a certain portion of the
dissertation \cite{Khan_Dissertation}. Also, the referee of
\cite{Khan_Smoothable4} deserves thanks for suggesting Mostow
rigidity for the fibers $S$. Finally, the author appreciates his
discussions with Tom Farrell and Bruce Hughes, from which the
relevant techniques for the fibering result were made clear.

\bibliographystyle{alpha}
\bibliography{FiberingSplitting5}

\end{document}